\documentclass{article}

\PassOptionsToPackage{%
	style	= alphabetic%
}{biblatex}
\usepackage[doi=false,isbn=false,url=false,
	hyperref=auto,
	sorting=nyt,
	maxnames=10,
	maxcitenames=3,
	backend=biber,
	texencoding=auto,
	giveninits=true,
	block=space,
]{biblatex}

\DefineBibliographyStrings{english}{%
	andothers = {et al}
} 
 
\DeclareFieldFormat{eprint:arxiv}{%
	\em \href{https://arxiv.org/abs/#1}{arXiv:\allowbreak #1}}

\DeclareFieldFormat{eprint:mrnumber}{%
	\href{http://www.ams.org/mathscinet-getitem?mr=MR#1}{MR\allowbreak #1}}

\DeclareFieldFormat{eprint:jstor}{%
	\href{http://www.jstor.org/stable/#1}{JSTOR\allowbreak #1}}

\DeclareFieldFormat{eprint:customeprint}{%
	\em Available at \upshape \ttfamily \href{https://#1}{#1}}

\DeclareFieldFormat{eprint:inprep}{%
	\em In preparation}

\DeclareFieldFormat{eprint:note}{%
	\em {#1}}

\DeclareFieldFormat{eprint:onarxiv}{%
	\em Available on arXiv}

\DeclareFieldFormat{eprint:toappear}{%
	\mbox{\em To appear}}

\DeclareSourcemap{
	\maps[datatype=bibtex]{
		\map{
			\step[fieldsource=mrnumber,		fieldtarget=eprint, final]
			\step[fieldset=eprinttype,		fieldvalue=mrnumber]
		}
		\map{
			\step[fieldsource=arxiv,		fieldtarget=eprint, final]
			\step[fieldset=eprinttype,		fieldvalue=arxiv]
		}
		\map{
			\step[fieldsource=jstor,		fieldtarget=eprint, final]
			\step[fieldset=eprinttype,		fieldvalue=jstor]
		}
		\map{
			\step[fieldsource=customeprint,	fieldtarget=eprint, final]
			\step[fieldset=eprinttype,		fieldvalue=customeprint]
		}
		\map{
			\step[fieldsource=inprep,		fieldtarget=eprint, final]
			\step[fieldset=eprinttype,		fieldvalue=inprep]
		}
		\map{
			\step[fieldsource=note,			fieldtarget=eprint, final]
			\step[fieldset=eprinttype,		fieldvalue=note]
		}
		\map{
			\step[fieldsource=onarxiv,		fieldtarget=eprint, final]
			\step[fieldset=eprinttype,		fieldvalue=onarxiv]
		}
		\map{
			\step[fieldsource=toappear,		fieldtarget=eprint, final]
			\step[fieldset=eprinttype,		fieldvalue=toappear]
		}
	}
}

\DeclareFieldFormat{doi}{%
	\href{https://doi.org/#1}{DOI}%
}

\DeclareFieldFormat{url}{%
	\href{https://#1}{\texttt{#1}}%
}

\DeclareFieldFormat{comments}{%
	\em #1%
}

\DeclareFieldFormat
	[article,inbook,incollection,inproceedings,patent,thesis,unpublished]
	{title}{{#1\isdot}}

\DeclareFieldFormat
	[article,book,inbook,incollection,inproceedings,patent,thesis,unpublished, online]
	{date}{\mkbibparens{#1}}
\DeclareFieldFormat
	[article,book,inbook,incollection,inproceedings,patent,thesis,unpublished]
	{volume}{\mkbibbold{#1}}
\DeclareFieldFormat
	{pages}{\mkbibparens{#1}}

\DeclareFieldFormat{edition}{%
	\ifinteger{#1}
	{\mkbibordedition{#1}~\bibstring{edition}\addcomma}
	{#1~\bibstring{edition}\addcomma}}

\renewbibmacro*{volume+number+eid}{%
	\printfield{volume}%
\setunit*{\adddot}%
	\printfield{number}%
\setunit{\space}%
	\printfield{eid}%
}

\DeclareBibliographyDriver{article}{%
	\usebibmacro{bibindex}%
	\usebibmacro{begentry}%
	\usebibmacro{author}%
\setunit{\addspace}%
	\usebibmacro{date}%
\newunit\newblock
	\usebibmacro{title}%
\newunit\newblock
	\printfield{journaltitle}\addperiod%
\setunit*{\addspace}%
	\usebibmacro{volume+number+eid}%
\setunit*{\addspace}\newblock
	\printfield{pages}
\setunit*{\addspace}
	\printfield{comments}%
	\usebibmacro{eprint}%
\setunit*{\addnbspace}%
	\printfield{doi}%
	\usebibmacro{finentry}%
}

\DeclareBibliographyDriver{online}{%
	\usebibmacro{bibindex}%
	\usebibmacro{begentry}%
	\usebibmacro{author}%
\setunit{\addspace}%
	\usebibmacro{date}%
\newunit\newblock
	\usebibmacro{title}%
\newunit\newblock
	\usebibmacro{eprint}%
	\usebibmacro{finentry}%
}

\DeclareBibliographyDriver{book}{%
	\usebibmacro{bibindex}%
	\usebibmacro{begentry}%
	\usebibmacro{author}%
\setunit{\addspace}\newblock
	\usebibmacro{date}%
\newunit\newblock
	\usebibmacro{title}%
\setunit*{\addspace}\newblock
	\printfield{edition}%
\newunit\newblock
	\printlist{publisher}
\setunit*{\addspace}%
	\printfield{volume}%
\setunit*{\addspace}\newblock%
	\printfield{comments}%
	\usebibmacro{eprint}%
\setunit*{\addspace}
	\printfield{doi}
	\usebibmacro{finentry}%
}

\DeclareBibliographyDriver{thesis}{%
	\usebibmacro{bibindex}%
	\usebibmacro{begentry}%
	\usebibmacro{author}%
\setunit{\addspace}\newblock
	\usebibmacro{date}%
\newunit\newblock
	\usebibmacro{title}.%
\setunit*{\addspace}\newblock
	\space Thesis, \printlist{institution}
	\usebibmacro{eprint}%
	\printfield{comments}%
\setunit*{\addspace}
	\printfield{doi}%
	\usebibmacro{finentry}%
}

\DeclareBibliographyDriver{inproceedings}{%
	\usebibmacro{bibindex}%
	\usebibmacro{begentry}%
	\usebibmacro{author}%
\setunit{\addspace}%
	\usebibmacro{date}%
\newunit\newblock
	\usebibmacro{title}%
\newunit\newblock
	\printfield{booktitle}\addcomma%
\setunit*{\addspace}%
	\printfield{series}\addcomma%
\setunit*{\addspace}\newblock
	\usebibmacro{editor}\addcomma
\setunit*{\addspace}\newblock
	\printlist{publisher}
\setunit*{\addspace}%
	\printfield{volume}%
\setunit*{\addspace}%
	\printfield{pages}
\setunit*{\addspace}
	\usebibmacro{eprint}%
	\printfield{comments}%
\setunit*{\addnnspace}
	\printfield{doi}
	\usebibmacro{finentry}%
}

\DeclareBibliographyDriver{inbook}{%
	\usebibmacro{bibindex}%
	\usebibmacro{begentry}%
	\usebibmacro{author}%
\setunit{\addspace}%
	\usebibmacro{date}%
\newunit\newblock
	\usebibmacro{title}%
\newunit\newblock
	\printfield{booktitle}\addcomma%
\setunit*{\addspace}%
	\printfield{series}\addcomma%
\setunit*{\addspace}\newblock
	\usebibmacro{editor}\addcomma
\setunit*{\addspace}\newblock
	\printlist{publisher}
\setunit*{\addspace}%
	\printfield{volume}%
\setunit*{\addspace}%
	\printfield{pages}
\setunit*{\addspace}
	\usebibmacro{eprint}%
	\printfield{comments}%
\setunit*{\addspace}
	\printfield{doi}
	\usebibmacro{finentry}%
}

\DeclareBibliographyDriver{incollection}{%
	\usebibmacro{bibindex}%
	\usebibmacro{begentry}%
	\usebibmacro{author}%
\setunit{\addspace}%
	\usebibmacro{date}%
\newunit\newblock
	\usebibmacro{title}%
\newunit\newblock
	\printfield{booktitle}\addcomma%
\setunit*{\addspace}%
	\printfield{series}\addcomma%
\setunit*{\addspace}\newblock
	\printlist{publisher}
\setunit*{\addspace}%
	\printfield{volume}%
\setunit*{\addspace}%
	\printfield{pages}
\setunit*{\addspace}
	\usebibmacro{eprint}%
	\printfield{comments}%
\setunit*{\addspace}
	\printfield{doi}
	\usebibmacro{finentry}%
}

\AtEveryBibitem{%
	\clearfield{day}%
	\clearfield{month}%
} 
 
\DeclareNameAlias{sortname}{given-family}
\DeclareNameAlias{default}{given-family}

\usepackage{geometry}
\usepackage[hidelinks]{hyperref}
\hypersetup{
	colorlinks,
	linkcolor={blue},
	citecolor={ForestGreen},
	urlcolor={blue}
}

\usepackage{graphicx}
\usepackage[font = normal, labelfont = sf]{subcaption}
\graphicspath{{LPBL_Figures/}}
\usepackage[font = normal, labelfont = sf, labelsep = period, singlelinecheck = false]{caption}

\usepackage[UKenglish]{babel}
\usepackage{amsmath, amsthm, amssymb, mathrsfs, bm, mathtools}
\usepackage{wasysym}
\allowdisplaybreaks

\usepackage{titlesec}
\usepackage{xifthen}
\usepackage{xspace}

\usepackage{xparse}

\usepackage{titlesec}
\titleformat{\title}
	{\sffamily \Huge \bfseries \scshape \boldmath}{}{}{}
\titleformat{\section}
	{\sffamily \Large \bfseries \scshape \boldmath}{\thesection}{1em}{}
\titleformat{\subsection}
	{\sffamily \large \bfseries \scshape \boldmath}{\thesubsection}{1em}{}
\titleformat{\subsubsection}
	{\sffamily \normalsize \bfseries \scshape \boldmath}{\thesubsubsection}{1em}{}
\titleformat{\paragraph}[runin]
	{\sffamily \normalsize \bfseries \scshape \boldmath}{\theparagraph}{1em}{}
\titleformat{\subparagraph}[runin]
	{\sffamily \normalsize \bfseries \scshape \boldmath}{\thesubparagraph}{1em}{}

\usepackage{enumitem}

\setlist[enumerate]{%
	topsep	= \medskipamount,
	itemsep = 0pt,
	label	= \ensuremath{(\alph*)}
}
\setlist[itemize]{%
	topsep	= \medskipamount,
	itemsep	= 0pt,
	label	= \bcdot
}
\setlist[description]{%
	topsep	= \smallskipamount,		
	itemsep	= \smallskipamount,		
	font	= {\mdseries\itshape},	
}

\newcommand{\negphantom}[1]{\ifmmode\settowidth{\dimen0}{$#1$}\else\settowidth{\dimen0}{#1}\fi\hspace*{-\dimen0}}

\newenvironment{center-small}
	{\par\centering\medskip}
	{\par\medskip\noindent}

\usepackage[dvipsnames,hyperref]{xcolor}

\newcommand{\Quad}[1]{
	\mathchoice
	{\quad\text{#1}\quad}
	{\text{ #1 }}
	{\text{ #1 }}
	{\text{ #1 }}
}

\newcommand{\Qand}{\Quad{and}}
\newcommand{\Qas}{\Quad{as}}

\newcommand{\Qfor}{\Quad{for}}
\newcommand{\Qforall}{\Quad{for all}}

\newcommand{\Qif}{\Quad{if}}

\newcommand{\Qsince}{\Quad{since}}
\newcommand{\Qwhere}{\Quad{where}}
\newcommand{\Qwhen}{\Quad{when}}
\newcommand{\Qwith}{\Quad{with}}

\NewDocumentCommand{\tv}{sm}{%
	\IfBooleanT{#1}{\bigl}\lVert #2 \IfBooleanT{#1}{\bigr}\rVert_\mathsf{TV}%
}

\newcommand{\one}  [1]{\bm1\{ #1 \}}

\newcommand{\bcdot}{\ensuremath{\bm{\cdot}}}
\newcommand{\cq}{\coloneqq}

\newcommand{\abs}  [1]{| #1 |}
\newcommand{\absb} [1]{\bigl| #1 \bigr|}

\newcommand{\rbr} [1]{ ( #1 ) }
\newcommand{\rbb} [1]{\bigl( #1 \bigr)}

\newcommand{\rbbb}[1]{\biggl( #1 \biggr)}

\newcommand{\sbr} [1]{ [ #1 ] }
\newcommand{\sbb} [1]{\bigl[ #1 \bigr]}

\newcommand{\bra} [1]{ \{ #1 \} }
\newcommand{\brb} [1]{\bigl\{ #1 \bigr\}}

\newcommand{\floor}  [1]{\lfloor #1 \rfloor}

\newcommand{\ceil}[1]  {\lceil #1 \rceil}

\newcommand{\midb}{\bigm|}

\newcommand{\Oh}  [1]{\mathcal{O}( #1 )}

\newcommand{\oh}  [1]{o( #1 )}

\DeclarePairedDelimiterXPP{\PR}[2]
	{\mathbb P_{#1}}{(}{)}{}%
	{#2}
\NewDocumentCommand{\pr}{som}{%
	\IfBooleanTF{#1}
	{\PR[\big]{\IfValueT{#2}{#2}}{#3}}
	{\PR{\IfValueT{#2}{#2}}{#3}}
}

\DeclarePairedDelimiterXPP{\EX}[2]
	{\mathbb E_{#1}}{[}{]}{}%
	{#2}
\NewDocumentCommand{\ex}{som}{%
	\IfBooleanTF{#1}
	{\EX[\big]{\IfValueT{#2}{#2}}{#3}}
	{\EX{\IfValueT{#2}{#2}}{#3}}
}

\DeclarePairedDelimiterXPP{\VAR}[2]
	{\mathbb V\mathrm{ar}_{#1}}{(}{)}{}%
	{#2}
\NewDocumentCommand{\var}{som}{%
	\IfBooleanTF{#1}
	{\VAR[\big]{\IfValueT{#2}{#2}}{#3}}
	{\VAR{\IfValueT{#2}{#2}}{#3}}
}

\let\originalexp\exp
\let\exp\relax
\DeclarePairedDelimiterXPP{\EXP}[1]{\originalexp}{(}{)}{}{#1}
\NewDocumentCommand{\exp}{sm}{%
	\IfBooleanTF{#1}
		{\EXP[\big]{#2}}
		{\EXP{#2}}
}

\DeclarePairedDelimiterX\SET[1]{\{}{\}}%
	{#1}
\NewDocumentCommand{\set}{sm}{%
	\IfBooleanTF{#1}
		{\SET[\big]{#2}}
		{\SET{#2}}
}

\usepackage{calc}
\newlength{\halfplusheight}
\NewDocumentCommand{\sumt}{smo}{%
	\mathchoice%
	{\textstyle\sum_{#2}\IfBooleanT{#1}{^\star}\IfValueT{#3}{^{#3}}\displaystyle}%
	{\sum_{#2}\IfBooleanT{#1}{^\star}\IfValueT{#3}{^{#3}}}%
	{\sum_{#2}\IfBooleanT{#1}{^\star}\IfValueT{#3}{^{#3}}}%
	{\sum_{#2}\IfBooleanT{#1}{^\star}\IfValueT{#3}{^{#3}}}%
}

\NewDocumentCommand{\sumd}{smo}{%
	\sum_{#2}\IfBooleanTF{#1}{^\star}{\IfValueT{#3}{^{#3}}}%
}

\newcommand{\intt}[2][]{
	\mathchoice
	{\ifthenelse{\isempty{#1}}
		{\textstyle \int_{#2}      \displaystyle}
		{\textstyle \int_{#2}^{#1} \displaystyle}}
	{\ifthenelse{\isempty{#1}}
		{\int_{#2}}
		{\int_{#2}^{#1}}}
	{\ifthenelse{\isempty{#1}}
		{\int_{#2}}
		{\int_{#2}^{#1}}}
	{\ifthenelse{\isempty{#1}}
		{\int_{#2}}
		{\int_{#2}^{#1}}}
}

\NewDocumentCommand{\prodt}{mo}{%
	\mathchoice%
	{\textstyle \prod_{#1}\IfValueT{#2}{^{#2}} \displaystyle}%
	{\prod_{#1}\IfValueT{#2}{^{#2}}}%
	{\prod_{#1}\IfValueT{#2}{^{#2}}}%
	{\prod_{#1}\IfValueT{#2}{^{#2}}}%
}

\NewDocumentCommand{\prodd}{smo}{%
	\prod_{#2}\IfBooleanT{#1}{^\star}\IfValueT{#3}{^{#3}}%
}

\DeclareMathOperator{\Exp}{Exp}

\DeclareMathOperator{\Pois}{Pois}

\DeclareMathOperator{\Bin}{Bin}

\DeclareMathOperator{\N}{N}

\newcommand{\RW}{\ifmmode \mathsf{RW} \else \textsf{RW}\xspace \fi}

\newcommand{\TV}{\ifmmode \mathsf{TV} \else \textsf{TV}\xspace \fi}

\newcommand{\tmix}{t_\mix}

\newcommand{\mix}{\textnormal{mix}}

\newcommand{\mbn}{\mathbb{N}}

\newcommand{\mbp}{\mathbb{P}}

\newcommand{\mbr}{\mathbb{R}}

\newcommand{\mbz}{\mathbb{Z}}

\newcommand{\mcf}{\mathcal{F}}

\newcommand{\mcs}{\mathcal{S}}

\newcommand{\mcx}{\mathcal{X}}

\newcommand{\mfS}{\mathfrak{S}}

\newcommand{\nt}{\addtocounter{equation}{1}\tag{\theequation}}

\newcommand{\eps}{\varepsilon}

\usepackage{manyfoot}
\SetFootnoteHook{\hspace*{-1.8em}}
\DeclareNewFootnote{bl}[gobble]
\newcommand{\blfootnote}[1]{\footnotebl{\sffamily#1}}

\def\IfAmpersandUseAlign#1#2&#3\EndIfAmpersandUseAlign%
	{%
		\if\relax\detokenize{#3}\relax
		\begin{equation*}%
		#1%
		\end{equation*}%
		\else
		\begin{align*}%
		#1%
		\end{align*}%
		\fi
	}
	\def\[#1\]%
	{%
		\IfAmpersandUseAlign{#1}#1&\EndIfAmpersandUseAlign
	}

\usepackage{eqparbox}

\newcommand{\qedtriangle}{\renewcommand{\qedsymbol}{\ensuremath{\triangle}}}

\usepackage[capitalise]{cleveref}

\crefname{figure}{Figure}{Figures}

\numberwithin{equation}{section}

\crefformat{equation}{\textup{(#2#1#3)}}
\crefmultiformat{equation}{\textup{(#2#1#3}}{\textup{, #2#1#3)}}{\textup{, #2#1#3}}{\textup{, #2#1#3)}}
\crefrangeformat{equation}{\textup{(#3#1#4--#5#2#6)}}

\crefformat{section}{\S#2#1#3}
\crefformat{subsection}{\S#2#1#3}
\crefformat{subsubsection}{\S#2#1#3}

\newenvironment{Proof}[1][\proofname]{%
	\proof[\upshape\bfseries\sffamily\boldmath#1]
}{\endproof}

\usepackage{etoolbox}
\usepackage{needspace}

\newtheoremstyle{sfsl}
	{1\baselineskip}		
	{1\baselineskip}		
	{\slshape}				
	{}						
	{\bfseries\sffamily}	
	{.}						
	{0.5em}					
	{\thmname{#1}\thmnumber{ #2}\thmnote{ {\mdseries(#3)}}}

\newtheoremstyle{sfup}
	{1\baselineskip}		
	{1\baselineskip}		
	{\upshape}				
	{}						
	{\bfseries\sffamily}	
	{.}						
	{0.5em}					
	{\thmname{#1}\thmnumber{ #2}\thmnote{ {\mdseries(#3)}}}

\theoremstyle{sfsl}

\newtheorem*{thm*}{Theorem}
\newtheorem{thm} {Theorem}[section]
\crefname{thm}{Theorem}{Theorems}

\newtheorem*{introthm*}{Theorem}
\newtheorem{introthm}{Theorem}

\crefname{introthm}{Theorem}{Theorems}

\newtheorem*{cor*}{Corollary}
\newtheorem{cor} [thm]{Corollary}
\crefname{cor}{Corollary}{Corollaries}

\newtheorem*{introcor*}{Corollary}

\crefname{introcor}{Corollary}{Corollaries}

\newtheorem*{introconj*}{Conjecture}

\crefname{introconj}{Conjecture}{Conjectures}

\newtheorem*{introques*}{Question}

\crefname{introques}{Question}{Questions}

\newtheorem*{lem*}    {Lemma}
\newtheorem{lem} [thm]{Lemma}
\crefname{lem}{Lemma}{Lemmas}

\newtheorem*{introlem*}{Lemma}

\crefname{introlem}{Lemma}{Lemmas}

\newtheorem*{prop*}    {Proposition}
\newtheorem{prop} [thm]{Proposition}
\crefname{prop}{Proposition}{Propositions}

\newtheorem*{clm*}    {Claim}

\crefname{clm}{Claim}{Claims}

\newtheorem*{defn*}    {Definition}
\newtheorem{defn} [thm]{Definition}
\crefname{defn}{Definition}{Definitions}

\newtheorem*{introdefn*}{Definition}
\newtheorem{introdefn}{Definition}

\crefname{introdefn}{Definition}{Definitions}

\newtheorem*{alg*}{Algorithm}

\crefname{alg}{Algorithm}{Algorithms}

\crefname{nota}{Notation}{Notations}

\newtheorem{innercustomgeneric}{\customgenericname}
\providecommand{\customgenericname}{}

\newcommand{\customgenname}{}

\newtheorem*{conj*}   {Conjecture}

\crefname{conj}{Conjecture}{Conjectures}

\newenvironment{conj-ind*}
	{\begin{quote}\textsf{\textbf{Conjecture.}}\slshape}
	{\end{quote}}
\newenvironment{conj-ind}
	{\begin{quote}\vspace{-\glueexpr\baselineskip+\topsep}\begin{customconj}}
	{\end{customconj}\end{quote}}

\newenvironment{question-ind*}
	{\begin{quote}\textsf{\textbf{Question.}}\slshape}
	{\end{quote}}
\newenvironment{question-ind}
	{\begin{quote}\vspace{-\glueexpr\baselineskip+\topsep}\begin{customquestion}}
	{\end{customquestion}\end{quote}}

\newenvironment{openproblem-ind*}
	{\begin{quote}\textsf{\textbf{Open Problem.}}\slshape}
	{\end{quote}}
\newenvironment{openproblem-ind}
	{\begin{quote}\vspace{-\glueexpr\baselineskip+\topsep}\begin{customopenproblem}}
	{\end{customopenproblem}\end{quote}}

\newtheorem*{hypothesis*}{Hypothesis}

\newtheorem*{hyp*}{Hypothesis}

\crefname{hyp}{Hypothesis}{Hypotheses}

\newtheorem*{rmk*}{Remark}

\theoremstyle{sfup}

\crefname{defn} {Definition}{Definitions}
\crefname{defnT}{Definition}{Definitions}

\newtheorem*{defnTT}{Definition}
\newenvironment{defnt*}
	{\pushQED{\qed}\renewcommand{\qedsymbol}{\ensuremath{\triangle}}\defnTT}
	{\popQED\enddefnTT}

\crefname{rmk} {Remark}{Remarks}
\crefname{rmkT}{Remark}{Remarks}

\newtheorem*{rmkTT}{Remark}
\newenvironment{rmkt*}
	{\pushQED{\qed}\renewcommand{\qedsymbol}{\ensuremath{\triangle}}\rmkTT}
	{\popQED\endrmkTT}

\crefname{rmks} {Remarks}{Remarks}
\crefname{rmksT}{Remarks}{Remarks}

\newtheorem*{rmks*} {Remarks}
\newtheorem*{rmksTT}{Remarks}
\newenvironment{rmkst*}
	{\pushQED{\qed}\renewcommand{\qedsymbol}{\ensuremath{\triangle}}\rmksTT}
	{\popQED\endrmksTT}

\newtheorem{intrormkT}{Remark}
	\newenvironment{intrormkt}
	{\pushQED{\qed}\renewcommand{\qedsymbol}{\ensuremath{\triangle}}\intrormkT}
	{\popQED\endintrormkT}

\crefname{intrormk} {Remark}{Remarks}
\crefname{intrormkT}{Remark}{Remarks}

\newtheorem*{intrormk*} {Remark}
\newtheorem*{intrormkTT}{Remark}
\newenvironment{intrormkt*}
	{\pushQED{\qed}\renewcommand{\qedsymbol}{\ensuremath{\triangle}}\intrormkTT}
	{\popQED\endintrormkTT}

\crefname{exm} {Example}{Examples}
\crefname{exmT}{Example}{Examples}

\newtheorem*{exm*} {Example}
\newtheorem*{exmTT}{Lemma}
	\newenvironment{exmt*}
	{\pushQED{\qed}\renewcommand{\qedsymbol}{\ensuremath{\triangle}}\exmTT}
	{\popQED\endexmTT}

\newtheorem*{note*} {Note}
\newtheorem*{noteTT}{Note}
	\newenvironment{notet*}
	{\pushQED{\qed}\renewcommand{\qedsymbol}{\ensuremath{\triangle}}\noteTT}
	{\popQED\endnoteTT}


\newcounter{parentnumber}

\makeatletter
\newenvironment{subtheorem-num}[1]{%
	\def\subtheoremcounter{#1}%
	\refstepcounter{#1}%
	\protected@edef\theparentnumber{\csname the#1\endcsname}%
	\setcounter{parentnumber}{\value{#1}}%
	\setcounter{#1}{0}%
	\expandafter\def\csname the#1\endcsname{\theparentnumber.\arabic{#1}}%
	\expandafter\def\csname theH#1\endcsname{thm.\theparentnumber.\arabic{#1}}%
	\unskip\ignorespaces
}{%
	\setcounter{\subtheoremcounter}{\value{parentnumber}}%
	\ignorespacesafterend
}
\makeatother

\makeatletter
\let\save@mathaccent\mathaccent
\newcommand*\if@single[3]{%
  \setbox0\hbox{${\mathaccent"0362{#1}}^H$}%
  \setbox2\hbox{${\mathaccent"0362{\kern0pt#1}}^H$}%
  \ifdim\ht0=\ht2 #3\else #2\fi
  }
\newcommand*\rel@kern[1]{\kern#1\dimexpr\macc@kerna}
\newcommand*\widebar[1]{\@ifnextchar^{{\wide@bar{#1}{0}}}{\wide@bar{#1}{1}}}
\newcommand*\wide@bar[2]{\if@single{#1}{\wide@bar@{#1}{#2}{1}}{\wide@bar@{#1}{#2}{2}}}
\newcommand*\wide@bar@[3]{%
  \begingroup
  \def\mathaccent##1##2{%
    \let\mathaccent\save@mathaccent
    \if#32 \let\macc@nucleus\first@char \fi
    \setbox\z@\hbox{$\macc@style{\macc@nucleus}_{}$}%
    \setbox\tw@\hbox{$\macc@style{\macc@nucleus}{}_{}$}%
    \dimen@\wd\tw@
    \advance\dimen@-\wd\z@
    \divide\dimen@ 3
    \@tempdima\wd\tw@
    \advance\@tempdima-\scriptspace
    \divide\@tempdima 10
    \advance\dimen@-\@tempdima
    \ifdim\dimen@>\z@ \dimen@0pt\fi
    \rel@kern{0.6}\kern-\dimen@
    \if#31
      \overline{\rel@kern{-0.6}\kern\dimen@\macc@nucleus\rel@kern{0.4}\kern\dimen@}%
      \advance\dimen@0.4\dimexpr\macc@kerna
      \let\final@kern#2%
      \ifdim\dimen@<\z@ \let\final@kern1\fi
      \if\final@kern1 \kern-\dimen@\fi
    \else
      \overline{\rel@kern{-0.6}\kern\dimen@#1}%
    \fi
  }%
  \macc@depth\@ne
  \let\math@bgroup\@empty \let\math@egroup\macc@set@skewchar
  \mathsurround\z@ \frozen@everymath{\mathgroup\macc@group\relax}%
  \macc@set@skewchar\relax
  \let\mathaccentV\macc@nested@a
  \if#31
    \macc@nested@a\relax111{#1}%
  \else
    \def\gobble@till@marker##1\endmarker{}%
    \futurelet\first@char\gobble@till@marker#1\endmarker
    \ifcat\noexpand\first@char A\else
      \def\first@char{}%
    \fi
    \macc@nested@a\relax111{\first@char}%
  \fi
  \endgroup
}
\makeatother

\usepackage{xparse}
\ExplSyntaxOn
\NewDocumentCommand{\mref}{m}{\quinn_mref:n {#1}}
\seq_new:N \l_quinn_mref_seq
\cs_new:Npn \quinn_mref:n #1
{
	\seq_set_split:Nnn \l_quinn_mref_seq { , } { #1 }
	\seq_pop_right:NN \l_quinn_mref_seq \l_tmpa_tl
	( 
	\seq_map_inline:Nn \l_quinn_mref_seq
	{ \ref{##1},\nobreakspace } 
	\exp_args:NV \ref \l_tmpa_tl 
	) 
}
\ExplSyntaxOff

\usepackage{stmaryrd}

\newcommand{\whp}{\ifmmode \mathsf{whp} \else \textsf{whp}\xspace \fi}

\DeclareMathOperator{\HG}{HG}

\newcommand{\ka}{\kappa}

\newcommand{\BL}{\ifmmode \mathsf{BL} \else \textsf{BL}\xspace \fi}
\newcommand{\RT}{\ifmmode \mathsf{RT} \else \textsf{RT}\xspace \fi}

\newcommand{\unif}[1][]{%
	\mathcal U%
	\ifthenelse{\equal{#1}{}}{}{_{#1}}%
}
\NewDocumentCommand{\Unif}{sm}{%
	\IfBooleanTF{#1}%
		{\mathcal U_{#2}}%
		{\operatorname{Unif}(#2)}%
}

\newcommand{\bb}[1]{\llbracket #1 \rrbracket}

\NewDocumentCommand{\DD}{smooo}{%
	d_{#2}%
	\IfValueT{#3}{%
		\IfBooleanT{#1}{\bigl}(%
		#3%
		\IfValueT{#4}{, \: #4}%
		\IfBooleanT{#1}{\bigr})%
	}%
}

\renewcommand{\P}{\mathbb{P}}

\usepackage{empheq}

\title{\sffamily%
	Limit Profile for the Bernoulli--Laplace Urn
}

\author{\sffamily Sam Olesker-Taylor\quad Dominik Schmid}
\date{\sffamily \today}

\begin{document}

\maketitle

\blfootnote{%
	Dominik Schmid%
\quad
	\href{mailto:d.schmid@uni-bonn.de}{d.schmid@uni-bonn.de}%
\hfill%
	\href{mailto:sam.olesker-taylor@warwick.ac.uk}{sam.olesker-taylor@warwick.ac.uk}%
\quad%
	Sam Olesker-Taylor%
\\
	Department of Mathematics, University of Bonn%
\hfill
	Department of Statistics, University of Warwick%
}

\vspace{-6ex}

\renewcommand{\abstractname}{\sffamily Abstract}

\begin{abstract}
\noindent
We analyse the convergence to equilibrium of the \textit{Bernoulli--Laplace urn model}:
\begin{itemize}[noitemsep, topsep=\smallskipamount]
	\item 
	initially,
	one urn contains $k$ red balls and a second $n-k$ blue balls;
	
	\item 
	in each step,
	a pair of balls is chosen uniform and their locations are switched.
\end{itemize}
Cutoff is known to occur at $\tfrac12 n \log \min\set{k, \sqrt n}$ with window order $n$ whenever $1 \ll k \le \tfrac12 n$.
We refine this by determining the \textit{limit profile}:
a function $\Phi$
such that
\begin{center-small}
\(\displaystyle
	d_\TV\rbb{ \tfrac12 n \log \min\set{k, \sqrt n} + \theta n }
\to
	\Phi(\theta)
\Qas
	n \to \infty
\Qforall
	\theta \in \mbr.
\)
\end{center-small}\noindent
Our main technical contribution, of independent interest, approximates a rescaled chain by a diffusion on $\mbr$ when $k \gg \sqrt n$, and uses its explicit law as a Gaussian process.
\end{abstract}


\small
\begin{quote}
\begin{description}
	\item [Keywords:]
	cutoff,
	limit profile,
	Bernoulli--Laplace,
	diffusion approximation
	
	\item [MSC 2020 subject classifications:]
	60J10, 60J60, 37A25, 37A30
	
\end{description}
\end{quote}
\normalsize







\sffamily
\boldmath
\setcounter{tocdepth}{1}
\tableofcontents
\unboldmath
\normalfont

\section{Introduction}

The focus of our work is on the famous \textit{Bernoulli--Laplace urn model}.
It consists of $n$ balls, $k$ of which are red and the remainder black, placed in two urns: $k$ in the first and the remainder in the second.
A single step of the process consist of choosing a pair of balls uniformly at random and swapping their positions.
It was introduced by Bernoulli and Laplace as a discrete model for the diffusion of two incompressible gases between two containers.

The process is a finite Markov chain, and converges to its (unique) equilibrium distribution.
The convergence is not smooth, but abrupt:
	it transitions from being far from mixed at time say $t_n - 100w_n$ to well mixed at $t_n + 100w_n$ in a \textit{window} of with order $w_n$ much smaller than the \textit{mixing time} $t_n$ (ie, $w_n / t_n \to 0$ as $n \to \infty$).
This phenomenon is called \textit{cutoff}, and was established by \textcite{DS:bernoulli-laplace} for the Bernoulli--Laplace model.
The model has received a great deal of attention subsequently,
	including \cite{DLS:bernoulli-laplace,S:many-urn-bl,CsST:gelfand-applications,LL:cutoff-ep-ip,FJ:spectrum-ep-ip,EN:bl-on-swaps,ABBHKS:bernoulli-laplace}, to give just a few.
We defer a detailed history to \S\ref{sec:intro:related-bl}.

We zoom into this window and determine the exact shape of the \textit{limit profile}.
We determine how far the process from equilibrium after time $t_n(\theta) \cq t_n + \theta w_n$ in the limit $n \to \infty$ as a function of $\theta \in \mbn$.
Establishing cutoff requires only analysis of the limits $\theta \to \pm\infty$.

Obtaining the limit profile for Markovian models which exhibit cutoff is a recent topic.
Usually, cutoff is established via upper and lower bounds using highly distinct methods. This does not lend itself well to the precision needed for the limit profile.
We mention particularly a recent breakthrough by \textcite{T:limit-profile} establishing the limit profile for the random-transposition card shuffle via a new approximation lemma, closing a long-standing open question.
We defer further discussion, including on available techniques, to~\S\ref{sec:intro:previous-lp} and~\S\ref{sec:intro:methodsBL}.

As we explain in \S\ref{sec:intro:methodsBL}, we could not make any of these approaches work.
Instead we take a completely different and, to the best of our knowledge, new approach to limit profiles here. It is based on scaling limits for the number of red balls in the first urn, which depends on $\lim k/\sqrt n$.
For $k \gg \sqrt n$, we show (weak) convergence of the process to an Ornstein--Uhlenbeck diffusion using tools from stochastic analysis,	and take advantage of the explicit law of this Gaussian process.
For $k \asymp \sqrt n$, we approximate the birth--death chain by an $M/M/ \infty$ queue.
For $k \ll \sqrt n$, we relate the analysis to the coupon-collector problem.

\subsection{Model and Main Result}

We now define the model and limit profiles precisely and state the main result.

\begin{introdefn}[Bernoulli--Laplace Urn]
Let $n, k \in \mbn$ with $k < n$.
The \textit{Bernoulli--Laplace urn model}
consists
of two urns and $n$ balls, $k$ of which are red and $n-k$ black;
initially, $k$ balls, chosen arbitrarily, are placed in the first urn and the remainder in the second.
At the incidents of a rate-1 Poisson process,
pick two balls uniformly at random:
	if they are in the same urn, do nothing;
	otherwise, switch their positions.

Denote the number of red balls in the first urn after time $t \ge 0$ by $X_t^n$.
	%
\end{introdefn}



This definition is change equivalent, up to a time change, to picking at rate $1$ a ball from each urn independently, and swapping them. 
The time-$t$ state of the process is completely described by the number $X_t^n$ of red balls in the first urn.
Indeed, there are then $k - X_t^n$ black balls in the first urn, $k - X_t^n$ red in the second and $n - k + X_t^n$ black in the second.

It is intuitively clear that the balls become uniformly distributed over the two urns as time progresses, subject to the constraint that the first urn always contains $k$ balls and the second $n-k$.
In other words, the equilibrium distribution $\pi_n$ of $(X_t^n)_{t\ge0}$ corresponds to uniform sampling $k$ times from $n$ balls without~replacement:
\[
	\pi_n(r)
=
	\textstyle
	\binom{k}{r} \binom{n-k}{n-r} \big/ \binom{n}{k}
=
	\pr*{ \HG(n, k, k) = r },
\]
where $\HG(n, k, k)$ is the \textit{hypergeometric} distribution with size $n$ population, $k$ success states and $k$ draws.
It is natural to ask how many swaps are needed to be close to equilibrium.

This question of the \textit{mixing time} was first studied by \textcite{DS:bernoulli-laplace}.
They showed that time $t_\theta \cq \tfrac14 n (\log n + \theta)$  is necessary and sufficient in the limit $n \to \infty$, assuming $k$ to be of order $n$.
In other words, the law of $X_{t_\theta}^n$ is far from $\pi_n$ if $\theta$ is large and \emph{negative}, whilst it is close if $\theta$ is large and \emph{positive}.
This sharp transition is called \textit{cutoff}.

Our contribution is in determining the fluctuations around the mixing time $t_0$:
	we evaluate precisely the total-variation distance $d_\TV(t)$ between $X_t^n$ and $\pi_n$ in the limit as $n \to \infty$ when $t = t_\theta$ \emph{for every fixed $\theta \in \mbr$};
	this is the \textit{limit profile}.
Contrastingly, for \emph{cutoff} it suffices to show $d_\TV(t_\theta) \to 0$ as $\theta \to -\infty$ and $d_\TV(t_\theta) \to 1$ as $\theta \to +\infty$.
We make this precise now.

\begin{introdefn*}[Mixing Time]
Let $\Omega$ be a finite set and $\mu$ and $\pi$ be distributions on $\Omega$.
Then,
\[
	\tv{ \mu - \pi }
\cq
	\sup_{A \subseteq \Omega}
	\abs{ \mu(A) - \pi(A) }
\in
	[0,1]
\]
is the \textit{total-variation} (\textit{\TV}) distance between $\mu$ and $\pi$.
Let $(X_t)_{t\ge0}$ denote a continuous-time irreducible Markov chain on a finite state space $\Omega$ with equilibrium distribution $\pi$.
Then,
for $t \ge 0$ and $\eps \in (0,1)$, respectively,
\[
	d_\TV(t)
\cq
	\max_{x \in \Omega}
	\tv{ \pr{X_t \in \cdot \, | \, X_0=x } -  \pi }
\Qand
	\tmix(\eps)
\cq
	\inf\bra{ t \ge 0 \mid d_\TV(t) \le \eps }
\]
are the \textit{(worst-case) \TV distance to equilibrium} and \textit{(precision-$\eps$) mixing time}, respectively.
\end{introdefn*}

The convergence theorem for Markov chains \cite[Theorem~1.8.3]{N:markov} ensures that $\tmix(\eps) < \infty$ for all $\eps \in (0,1)$  whenever the Markov chain is finite and irreducible.

\begin{introdefn*}[Cutoff and Limit Profile]
Let $(X^n)_{n\in\mbn}$ be a sequence of irreducible finite-state-space Markov chains, and use a superscript-$n$ to indicate statistics for the $n$-th chain.
Then, the sequence exhibits \textit{cutoff} at time $(t_0^n)_{n\in\mbn}$ with \textit{window} at most $(w^n)_{n\in\mbn}$ if
\[
	\lim_{\theta\to\infty}
	\lim_{n\to\infty}
	d_\TV^n(t_\theta^n)
=
	0
\Qand
	\lim_{\theta\to-\infty}
	\lim_{n\to\infty}
	d_\TV^n(t_\theta^n) = 1
\Qwhere
	t^n_\theta
\cq
	t_0^n + \theta w^n
\Qfor
	\theta \in \mbr.
\]
Moreover, the function $\varphi : \mbr \to (0,1)$ is the \textit{limit profile} with respect to $(t_\theta^n)_{n\in\mbn}$ if
\[
	\varphi(\theta)
\cq
	\lim_{n\to\infty}
	d_\TV^n(t_\theta^n)
\quad
	\text{exists for all $\theta \in \mbr$ and takes values in $(0,1)$}.
\]
	%
\end{introdefn*}

Our main result is the evaluation of the limit profile for the Bernoulli--Laplace urn model.

\begin{introthm}[Limit Profile for Bernoulli--Laplace]
\label{res:intro:bl}
For each $n \in \mbn$, let $k_n \in \bra{1, ..., \floor{\tfrac12 n}}$.
Let $d_\TV^n(t)$ denote the \TV to equilibrium for the Bernoulli--Laplace urn with parameters $(k_n, n)$ after time $t$.
We distinguish three cases according to the asymptotic behaviour of $k_n / \sqrt n$.
\begin{enumerate}
\item 
Suppose that $k_n / \sqrt n \to \infty$ as $n \to \infty$.
Then,
\[
	d^n_\TV\rbb{ \tfrac14 n \log n + \theta n }
\to
	\tv{ N(e^{-2\theta}, 1) - N(0, 1) }
\Qas
	n \to \infty,
\]
where $N$ denotes the law of a Gaussian random variable.  

\item 
Suppose that $k_n / \sqrt n \to \sqrt \alpha$ for some $\alpha>0$ as $n \to \infty$.
Then,
\[
	d_\TV\rbr{ \tfrac14 n \log n + \theta n }
&\to
	\tv{ \Pois\rbr{ \alpha + \sqrt \alpha e^{-2\theta} } - \Pois(\alpha) }
\Qas
	n \to \infty.
\intertext{where $\Pois$ denotes the law of a Poisson random variable.
Reparametrising,}
	d^n_\TV\rbb{ \tfrac12 n \log k_n + \theta n }
&\to
	\tv{ \Pois(\alpha + e^{-2\theta}) - \Pois(\alpha) }
\Qas
	n \to \infty.
\]

\item 
Suppose that $k_n / \sqrt n \to 0$ and $k_n \to \infty$ as $n \to \infty$.
Then,
\[
	d^n_\TV\rbb{ \tfrac12 n \log k_n + \theta n }
\to
	\pr{ G(0,1) > 2\theta }
\Qas
	n \to \infty,
\]
where $G \sim \operatorname{Gum}(0,1)$ is a standard Gumbel random variable.
\end{enumerate}
\end{introthm}

\begin{intrormkt}
Let us verify that the claimed limit in middle case ($k_n \asymp \sqrt n$) approximates the outer cases when $\alpha$ is very large or small.

We start with $\alpha$ large,
and use the normal approximation to the Poisson:
\[
	\alpha^{-1/2} \rbb{ \Pois(\alpha + \sqrt \alpha e^{-2\theta}) - \alpha }
\to
	N(e^{-2\theta}, 1)
\Qand
	\alpha^{-1/2} \rbb{ \Pois(\alpha) - \alpha }
\to
	N(0, 1).
\]
noting that $(\alpha + \sqrt \alpha e^{-2\theta})/\alpha \to 1$ as $\alpha \to \infty$.
Hence,
\[
	\lim_{\alpha\to\infty}
	\tv{ \Pois(\alpha + \sqrt \alpha e^{-2\theta}) - \Pois(\alpha) }
=
	\tv{ N(e^{-2\theta}, 1) - N(0, 1) }.
\]

We now turn to $\alpha$ small,
and use the Bernoulli approximation to the Poisson:
\[
&	\lim_{\alpha\to0}
	\tv{ \Pois(\alpha + e^{-2\theta}) - \Pois(\alpha) }
=
	\pr*{ \Pois(e^{-2\theta}) = 0 }
=
	e^{-2^{2\theta}}
=
	\pr*{ G > 2\theta }.
\qedhere
\]
\end{intrormkt}

%
%
%

\subsection{Related Work on the Bernoulli--Laplace model}
\label{sec:intro:related-bl}

The Bernoulli--Laplace model is analogous to the \textit{$k$-particle exclusion process} on the complete graph on $[n] \cq \set{1, ..., n}$.
This arises naturally as a projection of the \textit{interchange process}:
\begin{itemize}[noitemsep]
	\item 
	assign to each vertex $v \in [n]$ in the complete graph a unique particle labelled $\ell \in [n]$;
	
	\item 
	at rate $1$, pick an edge uniformly at random, and swap the particles at the endpoints.
\end{itemize}
Let $I^n_t(\ell) \in [n]$ denote the location of particle $\ell \in [n]$ after time $t$.
Then, the state of the \textit{interchange process} at time $t$ is $I^n_t \cq \rbr{ I^n_t(\ell) }_{\ell=1}^n$.
The \textit{$k$-particle exclusion process} is the projection from the \emph{vector} $I^n_t$ to the \emph{set} $S^{n,k}_t \cq \set{ I^n_t(\ell) }_{\ell=1}^k$ of locations of labels $1$ through~$k$.

To compare with Bernoulli--Laplace, place vertices $\set{1, ..., k}$ in the first urn and $\set{k+1, ..., n}$ in the second, and label the red balls $1$ through $k$ and the black $k+1$ through $n$.
By the symmetry of the complete graph, it is not important \emph{which} of the vertices $1$ through $k$ are occupied by the particles with labels $1$ through $k$.
Projecting gives Bernoulli--Laplace.

Cutoff for both the exclusion and interchange processes on the complete graph was established via probabilistic methods based on couplings by \textcite{LL:cutoff-ep-ip}.
The idea behind the coupling for the exclusion process is that it takes time $\tfrac14 n \log n + \Oh n$ to hit the `bulk', from which it takes a further time $\Oh n$ to couple two realisations.
We follow a similar probabilistic approach in the present paper. However, we require significantly refined approximations in order to determine the exact behaviour of the $\Oh n$ terms.

It is, in fact, possible to write down the entire spectrum for both the exclusion and interchange processes on the complete graph.
Intuitively, the non-zero eigenvalues can be related to subsets of particles.
A short, beautiful proof was obtained recently by \textcite{FJ:spectrum-ep-ip}.
The symmetry of the process (precisely, \textit{vertex transitivity}) means that the distance to equilibrium \emph{measured in $\ell_2$} (not \TV, which is equivalent to $\ell_1$) admits a representation depending only on the \emph{eigenvalues}, not the (complicated) \emph{eigenvectors}, from which, given the spectrum, it is straightforward to find the limit profile \emph{for $\ell_2$ distance}.

A related, and simpler, urn model is the \textit{Ehrenfest urn}.
In the simplest case, there are two urns and $n$ balls.
In a single step, a ball is chosen uniformly at random and moved to the other urn.
In the multi-urn version, the chosen ball is moved to a uniformly chosen urn.
Cutoff for this model was established by \citeauthor{CsST:harmonic-analysis-finite-groups}~\cite{CsST:gelfand-applications,CsST:harmonic-analysis-finite-groups}, and the limit profile determined by  Nestoridi and the first author \cite[Theorem~C]{NOt:limit-profiles:rev}, both in the multi-urn set-up.
A Normal profile, analogous to that proved here for Bernoulli--Laplace with $k = \tfrac12 n$, is established after $\tfrac12 n (\log n + \theta)$ steps.

Focussing on the two-urn model, each urn contains approximately half the balls in equilibrium.
Hence, when this is the case, if consecutive steps take from different urns, this is very similar to swapping a pair of balls, as in the Bernoulli--Laplace step. Moreover, the order of the steps is exchangeable.
Thus, in some sense, $2t$ steps of Ehrenfest urn closely approximates $t$ steps of Bernoulli--Laplace with $k = \tfrac12 n$, once each urn has approximately half the balls.
This condition is typically satisfied after $\omega$ steps whenever $\omega \gg n$---in particular, well before the mixing time, which is order $n \log n$.
It is thus not surprising that the mixing time of the Bernoulli--Laplace urn when $k = \tfrac12 n$ is simply half that of the Ehrenfest urn.

\subsection{Previous Work Establishing Limit Profiles}
\label{sec:intro:previous-lp}

Establishing cutoff is a long-standing objective in the study of Markov chains. 
Establishing the limit profile is a much more recent affair, and gaining significant traction.
We are only aware of a handful of classical examples, such as the random walk on the hypercube \cite{DGB:hypercube-profile} (for which a short argument was given recently in \cite[Theorem~5.1]{NOt:limit-profiles:rev-arxiv}) and the famous riffle shuffle \cite{BD:riffle-shuffle}.
On the other hand, there are a many references published within the last five years. The techniques used are quite diverse.

\textcite{T:limit-profile} recently introduced an approximation lemma using representation theory of the symmetric group.
This was extended in several directions by the first author and Nestoridi \cite{NOt:limit-profiles:rev,NOt:limit-profiles:proj}.
Representation-theoretic and orthogonal-polynomial techniques have been used for
	many card shuffles \cite{T:limit-profile,NOt:limit-profiles:rev,FTW:quantum-rt,NOt:limit-profiles:proj,N:limit-profiles:comp},
	the Ehrenfest urn and a Gibbs sampler \cite{NOt:limit-profiles:rev}
as well as
	multi-allelic Moran mutation--reproduction models \cite{C:moran-limit-profile}
Classically probabilistic approaches are used for
	random walks on Ramanujan and random Cayley graphs \cite{LP:ramanujan,HOt:rcg:abe:cutoff}
as well as
	repeated averages \cite{CDSZ:repeated-avg}, and the exclusion process on the circle \cite{L:exclusion-circle-profile}.
Very recently, tools from integrable probability have been used to analyse
	the asymmetric exclusion processes on the segment \cite{BN:asep-profile,HS:limit-asep}
and
	Metropolis biased card shuffling \cite{Z:cutoff-biased-shuffle}.

\subsection{Methods for Establishing Limit Profiles for BL and Related Models}
\label{sec:intro:methodsBL}

Cutoff was established for Bernoulli--Laplace by \textcite{DS:bernoulli-laplace} almost 40 years ago.
\citeauthor{CsST:harmonic-analysis-finite-groups}~\cite{CsST:gelfand-applications,CsST:harmonic-analysis-finite-groups} generalised the framework to \textit{Johnson schemes}, of which Bernoulli--Laplace is an example.
All these proofs
are algebraic.
They rely on representation theory of the symmetric group $\mfS_n$, and the particular expression of the model as a \textit{Gelfand pair},
a special case of homogeneous spaces.

Representation theory was also used by \textcite{DS:random-trans} to establish cutoff for the \textit{random-transposition} card shuffle:
	two cards are chosen uniformly at random from a deck, and their positions swapped;
	this corresponds to multiplying by a uniformly random transposition, hence the name.
Much more recently, \textcite[Lemma~2.1]{T:limit-profile} developed a novel approach refining this method to capture the limit profile, not only cutoff.
It involved evaluating a weighted sum over an arbitrarily large collection of \textit{characters}.

Bernoulli--Laplace can be seen as a projection of the random-transposition shuffle, where only the \textit{colour} of the card is considered.
Indeed, random transpositions is exactly the interchange process discussed earlier.
Nestoridi and the first author extended \citeauthor{T:limit-profile}'s approach to random walks on homogeneous spaces corresponding to Gelfand pairs in \cite[Lemma~C]{NOt:limit-profiles:rev}.
It requires evaluating an arbitrarily large sum over so-called \textit{spherical functions}.


Later, the same authors introduced an alternative, but related, approach for analysing projections of random walks on groups \cite{NOt:limit-profiles:proj}.
This involves lifting the projected walk from the homogeneous space (quotient) $X = G/K$ to the group $G$, where $K$ is a subgroup of $G$, with initial state uniform on $K$.
This lifts the Bernoulli--Laplace model to the random-transpositions shuffle on $\mfS_n$, started from the uniform distribution on $\mfS_k \times \mfS_{n-k} \subseteq \mfS_n$.
The $k$-particle interchange process corresponds starting uniformly on $\mfS_1 \times ... \times \mfS_1 \times \mfS_{n-k} \subseteq \mfS_n$.

In both cases, a certain sum over irreducible representations arises.
Combining calculations from \cite{T:limit-profile} with some refined information on occurrences of fixed points in uniform permutations, the limit profile for the $k$-particle interchange process was established. While this suggests an similar approach to establish the limit profile of the Bernoulli--Laplace model, we are unable to evaluate the sum in the case of exclusion (corresponding to Bernoulli--Laplace) as establishing cancellation of the summands is highly delicate.

Therefore, we take a different route in this paper, and leave open the question of achieving the same result via algebraic methods.


%

%

\subsection{Outline of Proof}

This paper is structured as follows.
\begin{itemize}
	\item [\S\ref{sec:bd}]
	We recall a standard reduction of the mixing time to that of a birth--death chains, and provide concentration results on the hitting time of set with a large stationary mass.
	
	\item [\S\ref{sec:gg}]
	This section is dedicated to the proof of the limit profile when $k/\sqrt n \to \infty$ as $n \to \infty$.
	We prove the convergence of the Bernoulli--Laplace urn to an Ornstein--Uhlenbeck process, which is a Gaussian process, from which we deduce the Gaussian limit profile.
	
	\item [\S\ref{sec:as}]
	We establish the limit profile when $k$ is of order $\sqrt n$. To do so, we identify an $M/M/\infty$ queue as the limit object of the rescaled  Bernoulli--Laplace urn.
	
	\item [\S\ref{sec:ll}]
	We conclude
	with the case where $k/\sqrt n \to 0$ as $n \to \infty$ by relating Bernoulli--Laplace, and its limit profile, to the coupon-collector problem. 
\end{itemize}

\subsection{Notation}

In the following, we collect some notation, which we frequently use throughout this paper.
For two functions $f,g : \mbn \to \mbr_+$,
with $f(x), g(x) \to \infty$ as $x \to \infty$,
we write
\begin{gather*}
	f \gg g
\Qif
	\lim_{x\to\infty}
	f(x)/g(x)
=
	\infty,
\qquad
	f \ll g
\Qif
	\lim_{x\to\infty}
	f(x)/g(x)
=
	0,
\\
	f \asymp g
\Qif
	\liminf_{x\to\infty}
	f(x)/g(x)
>
	0
\Qand
	\limsup_{x\to\infty}
	f(x)/g(x)
<
	\infty.
\end{gather*}
We write $f \gtrsim g$ whenever $f \gg g$ or $f \asymp g$. 
We use the Landon notation $\mathcal{O}$ and $o$ with respect to $n$, and we treat $k=k_n$ as a function of $n$. Moreover, it will be convenient to write
\[
	\ka_n \cq k_n / n
\Qfor
	n \in \mbn.
\]

\subsection{Acknowledgments}

We are especially grateful to Evita Nestoridi, with whom SOT had lengthy, detailed discussions around a representation-theoretic approach to the Bernoulli--Laplace model.
We also thank Lucas Teyssier, with whom SOT also had some similar discussions.

We thank Pedro Cardoso and Jimmy He for more general discussions in the area, including approximations by diffusions, and Konstantinos Panagiotou for pointing DS in the direction of the main diffusion-approximation theorem that we apply (namely, \cref{thm:Durrett}).

\section{Reduction to Birth--Death Chain}
\label{sec:bd}

We mentioned in the introduction that the number of red balls in the first urn characterises the process. Indeed, if there are $x$ red in the first, then there are $k_n - x$ black also in the first urn; the remaining balls are in the second.
The process counting the number of red balls in the first urn is a Markov chain taking values in $\bb{0, k_n} \cq \bra{0, 1, ..., k_n}$, and steps by $\pm1$ (or $0$). 
Exactly this reduction was used by \citeauthor{LL:cutoff-ep-ip}; see \cite[\S 2]{LL:cutoff-ep-ip} for more details. We write  $\P_x$ for the law of this Markov chain when starting from $x \in \bb{0,k_n}$.
\begin{lem}
\label{res:bd:rep}
The process $(X^n_t)_{t\ge0}$ counting the number of red balls in the first urn is~a birth--death chain.
Its birth- and death-rates are given by,
for~$x \in \bb{0, k_n}$,
\begin{alignat*}{2}
	q_n(x, x+1)
&=
	2 (k_n - x)^2 / n^2
&&\Quad{provided}
	x+1 \le k_n,
\\
	q_n(x, x-1)
&=
	2 x (n - 2 k_n + x) / n^2
&&\Quad{provided}
	x-1 \ge 0.
\end{alignat*}
In particular,
\[
	q_n(x, x+1) - q_n(x, x-1)
&=
	- 2 x / n,
\\
	q_n(x, x+1) + q_n(x, x-1)
&=
	4 \ka_n^2 (1-\ka_n)^2 + 8 (\tfrac12 - \ka_n)^2 x / n + 4 x^2 / n^2.
\]
\end{lem}

\begin{Proof}
Note that the number of red balls increases by $1$ if and only if we pick a black ball from the first urn and a red ball from the second urn. The probability of this event is exactly $q_n(x,x+1)$ when there are $x$ red balls in the first urn. Similarly, the number of red balls increases by $1$ if and only if we pick a black ball from the first urn and a red ball from the second urn. The probability of this event equals $q_n(x,x-1)$. Since we draw (with replacement) pairs of balls  at rate $1$, this gives the desired result. 
\end{Proof}


The behaviour of $(X^n_t)_{t\ge0}$ depends on the asymptotic behaviour of $k_n^2 / n = \ka_n^2 n$.
For example, the equilibrium distribution $\pi_n \sim \HG(n, k_n, k_n)$
has mean $k_n^2/n = \ka_n^2 n$ and variance
\[
	\frac{k_n^2}{n} \frac{n - k_n}{n} \frac{n - k_n}{n-1}
=
	\frac{\ka_n^2 (1 - \ka_n)^2 n}{1 - 1/n}
\sim
	\ka_n^2 n (1 - \ka_n)^2
\asymp
	\ka_n^2 n
=
	\frac{k_n^2}{n}.
\]

Let $C > 0$ be fixed.
If $k_n \ll \sqrt n$, then $\pi_n$ is concentrated on the singleton $\set{0}$. We are then interested in the time it takes to reach $C$. From there, the process is approximated by a coupon-collector.
Contrastingly, if $k_n \gg \sqrt n$, then the mean and variance diverge, and we approximate the behaviour in the by a suitable limiting object%
	---namely, an Ornstein--Uhlenbeck diffusion---%
once at $C$ standard deviations from the centre.
The times at which these points are hit in the two cases are, respectively,
\[
	T^-_n(C)
\cq
	\tfrac12 n \log k - \tfrac12 \log C
\Qand
	T^+_n(C)
\cq
	\tfrac14 n \log n - \tfrac12 \log C.
\]
The critical case $k_n \asymp \sqrt n$ can be analysed at either time, as
\(
	T^-_n(C k_n / \sqrt n)
=
	T^+_n(C).
\)
The limiting behaviour at this time is that of an $M/M/\infty$ queue.
We now formalise this intuition.

\begin{prop}
\label{res:bd:conc}
The following concentration holds.
\begin{itemize}
\item 
Suppose that $k_n \gtrsim \sqrt n$.
Then,
for all $\eps > 0$,
inspecting at time $T^+_n(C)$,
\[
	\limsup_{C\to\infty}
	\limsup_{n\to\infty}
	\pr*[k_n]{ 
		\abs{ 
			X^n_{T^+_n(C)}
		-	\rbr{ \ka_n^2 n + C \ka_n (1-\ka_n) \sqrt n }
		}
	>
		\eps \cdot C \ka_n \sqrt n
	}
=
	0.
\]
If $k_n \ll n$, then the same holds without the $1 - \ka_n$ factor in the shift.
If $k_n \asymp \sqrt n$,~then
\[
	\limsup_{C\to\infty}
	\limsup_{n\to\infty}
	\pr*[k_n]{ 
		\abs{ 
			X^n_{T^+_n(C)}
		-	C k_n / \sqrt n
		}
	>
		\eps \cdot C
	}
=
	0.
\]

\item 
Suppose that $k_n \lesssim \sqrt n$.
Then,
for all $\eps > 0$,
inspecting at time $T^-_n(C)$,
\[
	\limsup_{C\to\infty}
	\limsup_{n\to\infty}
	\pr*[k_n]{ 
		\abs{ X^n_{T^-_n(C)} - C }
	>
		\eps \cdot C
	}
=
	0.
\]
\end{itemize}
\end{prop}

To obtain this result, we need mean and variance estimates on $X^n_t$, rather than on $\pi_n$.
A simple inductive argument, given in \cite[Lemma~2.2]{LL:cutoff-ep-ip}, establishes exponential convergence of the mean $\ex[x]{X_t^n}$ to $k_n^2 / n = \ka_n^2 n$.
We repeat this short proof later, for completeness.

\begin{lem}
\label{res:bd:mean}
For all $t \ge 0$ and all $x \in \bb{0, k_n}$,
\[
	\ex[x]{X^n_t}
=
	\ka_n^2 n
+	(x - \ka_n^2 n) e^{-2t/n}.
\]
\end{lem}

\begin{lem}
\label{res:bd:var}
For all $t \ge 0$ and all $x \in \bb{0, k_n}$,
\[
	\var[x]{X^n_t}
&
=
	\ka_n^2 (1 - \ka_n)^2 n
	\rbr{1 - 1/n}^{-1}
	\rbb{ 1 - e^{-\frac4n(1-\frac1n)t} }
\\&
+	4 (\tfrac12 - \ka_n)^2
	(x - \ka_n^2 n) e^{-2t/n}
	\rbr{1 - 2/n}^{-1}
	\rbb{ 1 - e^{-\frac2n(1-\frac2n)t} }
\\&
+	\rbr{ (x - \ka_n^2 n) e^{-2t/n} }^2
	\rbb{ e^{4t/n^2} - 1 }.
\]
In particular,
if $t \le n (\log n)^2$, then
\[
	\var[\bar x]{\bar X^n_t}
\le
	\ka_n^2 n (2 + e^{-2t/n} n / k_n)
+	8 (k_n e^{-2t/n})^2 (\log n)^2 / n.
\]
\end{lem}

Given these lemmas, \cref{res:bd:conc} follows from Chebyshev's inequality.

\begin{Proof}[Proof of \cref{res:bd:conc}]
Suppose that $k_n \gtrsim \sqrt n$
(ie, $\ka_n \gtrsim 1/\sqrt n$).
By \cref{res:bd:mean,res:bd:var},
\[
	\ex[k_n]{X^n_t}
=
	\ka_n^2 n
+	C \ka_n (1 - \ka_n) \sqrt n
\Qand
	\var[k_n]{X^n_t}
\asymp
	\ka_n^2 n
	(1 + C)
\Qwhen
	t = T^+_n(C).
\]
We now apply Chebeyshev's inequality:
\[
	\pr*[k_n]{
		\abs{ \rbr{ X^n_{T^+_n(C)} - \ka_n^2 n } - C \ka_n (1-\ka_n) \sqrt n	}
	>
		\eps \cdot C \ka_n \sqrt n
	}
\lesssim
	\tfrac{1 + C}{\eps^2 C^2}
\to
	0
\Qas
	C \to \infty.
\]
If $k_n \ll n$ (ie, $\ka_n \ll 1$), then the same holds without the $1-\ka_n$ factor in the shift.

Now suppose that $k_n \lesssim \sqrt n$
(ie, $\ka_n \lesssim 1/\sqrt n$).
By \cref{res:bd:mean,res:bd:var},
\[
	\ex[k_n]{X^n_t}
=
	\ka_n^2 n
+	C (1 - \ka_n)
\Qand
	\var[k_n]{X^n_t}
\asymp
	1 + C
\Qwhen
	t = T^-_n(C).
\]
Again, we apply Chebyshev's inequality:
\[
	\pr*[k_n]{ 
		\abs{ X^n_{T^-_n(C)} - C }
	>
		\eps \cdot C
	}
\lesssim
	\tfrac{1+C}{\eps^2 C^2}
\to
	0
\Qas
	C \to \infty,
\]
using the fact that $\ka_n^2 n \lesssim 1$, so both $\ka_n^2 n$ and $C \ka_n$ are absorbed by $C \eps$ when $C \to \infty$.
	%
\end{Proof}

It remains to give the deferred proofs of \cref{res:bd:mean,res:bd:var}.
These rely on a standard result using the generator to evaluate the rate of change of the expected value of a statistic.

\begin{lem}
\label{res:bd:der-gen}
Let $Z = (Z_t)_{t\ge0}$ be an continuous-time Markov chain on a finite state space $\Omega$ with instantaneous transition-rates matrix $Q$.
Let $f : \Omega \to \mbr$.
Then,
\[
	\tfrac d{dt}
	\ex[z]{f(Z_t)}
=
	\ex[z]{ (Qf)(Z_t) }
\Qforall
	t \ge 0
\Qand
	z \in \Omega.
\]
\end{lem}

For convenience, we analyse the recentred version of $X$:
\[
	\bar X^n_t
\cq
	X^n_t - \ka_n^2 n.
\]
\cref{res:bd:mean} is proved by taking $f \cq f_1 : \bar x \mapsto \bar x$ and \cref{res:bd:var} by $f \cq f_2 : \bar x \mapsto \bar x^2$.

\begin{Proof}[Proof of \cref{res:bd:mean}]
Let $f_1 : \bar x \mapsto \bar x$.
Let $F_1(t) \cq \ex[\bar x]{f_1(X_t)} = \ex[\bar x]{X_t}$ for $t \ge 0$.
Then,
\[
	(Q f_1)(\bar x)
&
=
	q_n(\bar x, \bar x+1)
	\rbb{ f_1(\bar x+1) - f_1(\bar x) }
+	q_n(\bar x, \bar x-1)
	\rbb{ f_1(\bar x-1) - f_1(\bar x) }
\\&
=
	q_n(\bar x, \bar x+1) - q_n(\bar x, \bar x-1)
=
	- 2 \bar x / n,
\]
using \cref{res:bd:rep}.
Plugging this into \cref{res:bd:der-gen},
\[
	F_1'(t)
=
	- \tfrac2n
	F_1(t).
\]
Solving this differential equation, using $F_1(0) = \ex[\bar x]{\bar X^n_0} = \bar x$, proves the lemma.
\end{Proof}

\begin{Proof}[Proof of \cref{res:bd:var}]
Let $f_2 : \bar x \to \bar x^2$.
Let $F_2(t) \cq \ex[\bar x]{f_2(X_t)} = \ex[\bar x]{(\bar X^n_t)^2}$ for~\mbox{$t \ge 0$}.~%
Then,
\[
	(Q f_2)(\bar x)
&=
	( 2 \bar x + 1) q_n(\bar x, \bar x + 1)
+	(-2 \bar x + 1) q_n(\bar x, \bar x - 1)
\\&
=
	2 \bar x
	\rbb{ q_n(\bar x, \bar x + 1) - q_n(\bar x, \bar x - 1) }
+	\rbb{ q_n(\bar x, \bar x + 1) + q_n(\bar x, \bar x - 1) }
\\&
=
-	4 \bar x^2 / n
+	4 \bar x^2 / n^2 + 4 \ka_n^2 (1-\ka_n)^2 + 8 (\tfrac12 - \ka_n)^2 \bar x / n,
\]
using \cref{res:bd:rep}.
Plugging this into \cref{res:bd:der-gen},
\[
	F_2'(t)
&
=
-	\tfrac4n (1 - \tfrac1n)
	\ex[\bar x]{(\bar X^n_t)^2}
+	4 \ka_n^2 (1 - \ka_n)^2
+	\tfrac8n (\tfrac12 - \ka_n)^2
	\ex[\bar x]{X_t}
\\&
=
-	\tfrac4n (1 - \tfrac1n)
	F_2(t)
+	4 \ka_n^2 (1 - \ka_n)^2
+	\tfrac8n (\tfrac12 - \ka_n)^2
	\bar x e^{-2t/n},
\]
using \cref{res:bd:mean}.
We now use the integrating-factor method:
\[
	\tfrac d{dt}
	\rbb{ e^{\frac4n (1 - \frac1n) t} F_2(t) }
&
=
	e^{\frac4n (1 - \frac1n) t}
	\rbb{ F'_2(t) + \tfrac4n (1 - \tfrac1n) F_2(t) }
\\&
=
	4 \ka_n^2 (1-\ka_n)^2
	e^{\frac4n (1 - \frac1n) t}
+	\tfrac8n (\tfrac12 - \ka_n)^2 \bar x
	e^{\frac2n (1 - \frac2n) t}.
\]
Integrating,
and cancelling the $e^{\frac4n (1 - \frac1n) t}$ factor,
\[
	F_2(t)
&
=
	\frac{4 \ka_n^2 (1 - \ka_n)^2}{\frac4n (1 - \frac1n)}
+	\frac{\frac8n (\frac12 - \ka_n)^2 \bar x}{\frac2n (1 - \frac2n)}
	e^{-\frac2n t}
+	c
	e^{-\frac4n (1 - \frac1n) t}
\\&
=
	\frac{\ka_n^2 (1 - \ka_n)^2 n}{1 - 1/n}
+	\frac{4 (\frac12 - \ka_n)^2}{1 - 2/n}
	\bar x e^{-\frac2n t}
+	c
	e^{-\frac4n (1 - \frac1n) t},
\]
for some integration constant $c \in \mbr$ to be determined.
Namely, $F_2(0) = \bar x^2$, so
\[
	F_2(t)
=
	\frac{\ka_n^2 (1 - \ka_n)^2 n}{1 - 1/n}
	\rbb{ 1 - e^{-\frac4n (1 - \frac1n) t} }
+	\frac{4 (\frac12 - \ka_n)^2}{1 - 2/n}
	\bar x e^{-\frac2n t}
	\rbb{ 1 - e^{-\frac2n (1 - \frac1n) t} }
+	\bar x^2
	e^{-\frac4n (1 - \frac1n) t}.
\]
Now,
\(
	\var[\bar x]{\bar X^n_t}
=
	\ex[\bar x]{(\bar X^n_t)^2} - \ex[\bar x]{\bar X^n_t}^2
=
	F_2(t) - F_1(t)^2.
\)
So,
\[
	\var[\bar x]{\bar X^n_t}
&
=
	\ka_n^2 (1 - \ka_n)^2 n
	\rbb{ 1 - e^{-\frac4n (1 - \frac1n) t} }
/	(1 - 1/n)
\\&
+	4 (\tfrac12 - \ka_n)^2
	\bar x e^{-2t/n}
	\rbb{ 1 - e^{-\frac2n (1 - \frac1n) t} }
/	(1 - 2/n)
+	\rbr{ \bar x e^{-2t/n} }^2
	\rbb{ e^{4t/n^2} - 1 }.
\qedhere
\]

	%
\end{Proof}

\section{Limit Profile for $k \gg \sqrt n$}
\label{sec:gg}

Throughout this section, we assume that $\ka_n^2 n = k_n^2 / n \to \infty$ as $n \to \infty$.
In this regime, the equilibrium distribution is approximately $N(\ka_n^2 n, \ka_n^2 (1-\ka_n)^2 n)$.
We will rescale this.


\begin{Proof}[Outline of Proof]
When starting at a position which is in the standard fluctuations under the stationary distribution, we show that $(X^n_t)_{t \ge 0}$ is well-approximated by an Ornstein--Uhlenbeck process. We give the lower bound on the limit profile by applying the concentration result from \cref{res:bd:conc} on the birth--death chain reaching the start of the window, and then using the explicit description of the Ornstein--Uhlenbeck process as a Gaussian process. For the upper bound on the limit profile, we combine coupling and submartingale arguments to estimate the coalescence time of $(X^n_t)_{t \ge 0}$ with a stationary process.
\end{Proof}


\subsection{Convergence of the Bernoulli--Laplace Urn to a Diffusion}

In order to show convergence to a diffusion, we start by recalling the martingale problem introduced by Stroock and Varadhan; see \cite[\S 5.4]{KS:bm-sc} for an introduction.  

\begin{defn}[Martingale problem]
	Let $b,\sigma : \mbr \to \mbr$ be continuous functions.
	A stochastic process $(D_t)_{t \ge 0}$ solves the \textit{martingale problem} with respect to functions $b$ and $\sigma$ and initial condition $z \in \mbr$
	if 
	there is a unique probability measure $\mbp_z$, governing the law of the process $(D_t)_{t \ge 0}$, 
	such that $\pr[z]{D_0 = z} = 1$, and 
	\begin{equation}
		\label{eq:MartingaleProblem}
		\textup{d} D_t
	=
		b(D_t) \textup{d} t  + \sigma(D_t) \textup{d} B_t
	\end{equation}
	holds with respect to Brownian motion $(B_t)_{t \ge 0}$.
\end{defn}

A sufficient condition for the existence and uniqueness of the martingale problem is that $b$ and $\sigma$ are Lipschitz continuous, and we refer to \cite[\S 5.4]{KS:bm-sc} for a more detailed discussion. 
In the special case $\sigma(x) = A$ and $b(x) = -Bx$ for some constants $A,B>0$,
and all $x \in \mbr$, it is easy to check that 
the process $(D_t)_{t \ge 0}$ is an Ornstein--Uhlenbeck process with drift $B$ and diffusion constant $A$.
In particular,
if $D_0 = z$ almost surely,
then it is well-known that $D_t$ has an explicit Gaussian distribution:
\begin{equation}
\label{eq:OUExplicit}
	D_t
\sim 
	N\rbb{ z e^{-B t}, \tfrac{A^2}{2B} (1-e^{-2Bt}) }.
\end{equation}

The following convergence of a sequence of Markov chains to a diffusion can be found as a special case of \cite[Theorem~8.7.1]{D:stochastic-calc};
we omit the proof.

\begin{thm}
\label{thm:Durrett}
	For each $N \in \mbn$,
	let $(Z^N_t)_{t\ge0}$ be a continuous-time Markov chain on a finite state space $\mcs_N \subseteq \mbr$ with instantaneous transition-rates matrices $q_N = (q_N(x,y))_{x,y \in \mcs_N}$.
	Suppose there are Lipschitz continuous functions $b,\sigma : \mbr \to \mbr$ such that the following hold:
	\begin{itemize}
		\item 
		for any $(x_N)_{N \in \mbn}$ with $x_N \in \mcs_N$ and $x_N \to x$ for some $x \in \mbr$;
		
		\item 
		$\lim\limits_{N \to \infty}  \sum_{y : y \ne x_N} (y-x_N) q_N(x_N,y) \one{|x_N-y| \le 1} = b(x)$;
		
		\item 
		$\lim\limits_{N \to \infty} \sum_{y : y \ne x_N} (y-x_N)^2 q_N(x_N,y) \one{|x_N-y| \le 1} = \sigma(x)^2$.
	\end{itemize}
	Moreover, assume that,
	for all $R>0$ and all $\eps > 0$,
	\begin{equation}\label{eq:RadiusCondition}
		\lim_{N \to \infty}
		\sup_{|x| \le R}
		\sum_{y :|y - x| \ge \eps}
		q_N(x_N,y)
	=
		0. 
	\end{equation}
	Also assume that there are $z_N \in \mcs_N$ and $z \in \mbr$ such that
	\begin{equation}\label{eq:InitialStateCondition}
		Z^N_0 = z_N
	\Qforall
		N \in \mbn
	\Qand
		\lim_{N \to \infty}
		z_N
	=
		z .
	\end{equation}
	Then, for every fixed $T > 0$, in the sense of weak convergence, 
	\begin{equation}
	(Z^N_t)_{t \in [0,T]} \to (D_t)_{t \in [0,T]}, 
	\end{equation} where $(D_t)_{t \ge 0}$ solves the martingale problem with respect to $b$ and $\sigma$, with initial condition~$z$. 
\end{thm}

In the following, we apply \cref{thm:Durrett} when $k_n \gg \sqrt n$.
%
It will be convenient to rescale space by a factor of $\sqrt n \ka_n (1-\ka_n)$ and time by a factor of $n$ in order to ensure that we have convergence to a suitable diffusion. This is captured in the following lemma.

\begin{lem}\label{lem:ConvergenceSmaller}
	Fix $T>0$.
	Consider the birth--death chains $((X^n_t)_{t \ge 0})_{n \in \mbn}$ on $\bb{0, k_n}$ with initial state $X^n_0 = \ka_n^2 n+ z_n \sqrt n \ka_n(1-\ka_n)$ such that $z_n \to z$ for some constant $z>0$.
	Then,
	\begin{equation}
		\rbbb{ \frac{X^n_{nt} - \ka_n^2 n}{\sqrt n \ka_n (1-\ka_n)} }_{ t \in [0,T] }
	\to
		(\mathcal D_t)_{t \in [0,T] },
	\end{equation}
	where $(\mathcal D_t)_{t \in [0,T]}$ is the solution of the martingale problem in \cref{thm:Durrett} with $b(x)=-2x$ and $\sigma(x)=2$, and initial condition $z$. 
\end{lem} 

\begin{Proof}
Set $c_n \cq \sqrt n\ka_n(1-\ka_n) \in [\tfrac13 k_n/\sqrt n, k_n/\sqrt n]$ and consider the state space
\begin{equation}\label{def:SetSn}
	\mcs_n
=
	\bra{ x/c_n - \ka_n^2 n  \mid x \in \bb{0, k_n} }. 
\end{equation}
We denote by $(Y^n_t)_{t \ge 0}$ the continuous--time Markov chain on $\mcs_n$ given by
\begin{equation}
	Y^n_t
\cq
	\rbb{ X^n_{nt} - \ka_n^2 n } / c_n
\Qfor
	t \ge 0.
\end{equation}
We denote by $\widehat q_n = (\widehat q_n(y,y'))_{y, y' \in \mcs_n}$ the respective transition rates:
for all $y \in \mcs_n$
\begin{align*}
	\widehat q_n(y, y - c_n^{-1})
&=
	n q_n( c_n y + \ka_n^2 n, c_n y + \ka_n^2 n - 1 ),
\\
	\widehat q_n(y, y + c_n^{-1})
&=
	n q_n( c_n y + \ka_n^2 n, c_n y + \ka_n^2 n + 1 ),
\end{align*}
where the transition rates $q_n$ are from \cref{res:bd:rep}.
Let $y \in \mbr$.
Let $(x_n)_{n \in \N}$ be a sequence of initial conditions such that $x_n=\ka_n^2 n + c_n y_n \in \bb{0,k_n}$ for some $y_n \in \mbr$ with $y_n \to y$.
Then,
the following drift- and variance-limits hold courtesy of \cref{res:bd:rep}:
\[
	b(y)
&\cq
	\lim_{n \to \infty}
	\tfrac{1}{c_n}
	\rbb{ \widehat q_n(y_n,y_n + c_n^{-1}) - \widehat q_n(y_n,y_n - c_n^{-1}) }
= 
	- 2 y,
\\
	\sigma(y)^2
&\cq
	\lim_{n \to \infty}
	\tfrac{1}{c^2_n}
	\rbb{ \widehat q_n(y_n,y_n + c_n^{-1}) + \widehat q_n(y_n,y_n - c_n^{-1}) }
=
	4.
\]
The assumption \eqref{eq:RadiusCondition} on locality in \cref{thm:Durrett} is clearly satisfied for any fixed $R>0$ and $\eps > 0$. We conclude using \cref{thm:Durrett}. 
\end{Proof}
 
The stationary distribution $\pi_n$ is well approximated by a Gaussian distribution.
 
\begin{cor}\label{cor:StationaryDistribution}
	Let $\pi_n$ denote the stationary distribution of $X^n$ on $\bb{0,k_n}$.
	Let $X^n_\infty \sim \pi_n$ and $N \sim N(0,1)$ be a standard Gaussian random variable.
	Assume that $k_n \gg \sqrt n$.
	Then,
	\[
		\lim_{n \to \infty}
		\PR[\bigg]{}{ \frac{ X^n_\infty - \ka_n^2 n}{\sqrt n\ka_n(1-\ka_n)} \le x }
	= 
		\pr{N \le x}.
	\]
\end{cor} 

\begin{Proof}
	The result is immediate from Lemma~\ref{lem:ConvergenceSmaller} and  the observation \eqref{eq:OUExplicit} that the associated diffusion is an Ornstein--Uhlenbeck process, with $A = B = 2$, whose law convergences to a Gaussian with mean~$0$ and variance $1$ as $t \to \infty$. 
\end{Proof}

\subsection{Evaluating the Limit Profile for $k \gg \sqrt n$}
\label{sec:LimitLinear}

In order to establish the limit profile for $k_n \gg \sqrt n$, we first collect several basic observations on the birth--death chains $(X_t^n)_{t \ge 0}$. We consider the \textit{basic coupling} $\mathbf P$ between two birth--death chains $(X^n_t)_{t \ge 0}$ and $(\tilde X^n_t)_{t \ge 0}$ according to rates $q_n$ and initial conditions $x_0$ and $\tilde x_0$, respectively. In this coupling, we use the same Poisson clocks for both processes whenever the chains are in the same state. Otherwise, the two chains perform independent moves. We make the following observation on the ordering of birth--death chains. 

\begin{lem}\label{lem:StochasticDom}
	Let $(X^n_t)_{t \ge 0}$ and $(\tilde X^n_t)_{t \ge 0}$ be two birth--death chains evolving according to rates $q_n$ and with initial conditions $x_0^n$ and $\tilde x_0^n$, respectively, such that $x^n_0 \le \tilde x^n_0$. Then,
	\begin{equation*}
		X^n_t \le \tilde X^n_t
	\Qforall
		t \ge 0
	\Quad{under}
		\mathbf P.
	\end{equation*}
	Moreover, there exists some $C_0>0$ such that for all $C \ge C_0$, and all $n$ large enough, whenever $x^n_0 ,\tilde x^n_0 \in [\ka_n^2 n - C \ka_n \sqrt n ,  \ka_n^2 n + C \ka_n \sqrt n ]$, 
	\begin{equation*}
		\liminf_{n \to \infty}
		\mathbf P_{x^n_0, \tilde x^n_0}
		\rbb{ X^n_s ,\tilde X^n_s  \in [\ka_n^2 n - 2C \ka_n \sqrt n, \ka_n^2 n + 2C \ka_n \sqrt n ] \: \forall \, s \le n }
	\ge
		1 - 1/C.
	\end{equation*}
\end{lem}

\begin{Proof}
The first statement is immediate from the construction of the basic coupling $\mathbf P$. For the second statement, observe that the distance of $(X^n_t)_{t \ge 0}$ and $(\tilde X^n_t)_{t \ge 0}$ from $\ka_n^2 n$ is stochastically dominated by a simple random walk on $\mathbb{N}$ where the walker jumps at rate 
\[
	\max_{x \in [\ka^2_n n - 2C \ka_n \sqrt n, \ka^2_n n + 2C \ka_n \sqrt n]}
	\bra{ q_n(x,x-1)+ q_n(x,x+1) }
\in
	[\ka_n^2 , 4 \ka_n^2]
\]
for all $n$ large enough. The second statement is now immediate from a standard moderate deviation estimate for the symmetric simple random walk.
\end{Proof}


\begin{lem}\label{lem:Coalescence}
	Let $\delta>0$. Consider two birth--death chains $(X^n_t)_{t \ge 0}$ and $(\tilde X^n_t)_{t \ge 0}$ evolving according to rates $q_n$ and with initial conditions $x_0^n$ and $\tilde x_0^n$ such that
	\begin{equation}\label{eq:CloseToBoundary}
		\abs{ x^n_0 - \tilde x^n_0 }
	\le
		\delta \sqrt n \ka_n(1-\ka_n)
	\Qand
		x_0^n, \tilde x^n_0
	\in
		\sbb{ \ka^2_n n - \delta^{-1/4} \sqrt n\ka_n ,  \ka^2_n n + \delta^{-1/4} \sqrt n\ka_n }.
	\end{equation} 
	Then, there exists a constant $\delta_0>0$ such that if $\delta \in (0,\delta_0)$,
	then
	under the basic coupling~$\mathbf P$, 
	\begin{equation}\label{eq:GoalCoalFinal}
		\liminf_{n \to \infty}
		\mathbf P_{x^n_0, \tilde x^n_0}\rbb{ X^n_{4 \delta^{1/2} n} = \tilde X^n_{4 \delta^{1/2} n} }
	\ge
		1 - 2\delta^{1/4}. 
	\end{equation}
\end{lem}

\newcommand{\coal}{c}

\begin{Proof}
Let $\tau_\coal \cq \inf\bra{t \ge 0 \mid X^n_t = \tilde X^n_t}$.
The basic coupling $\mathbf P$ is \textit{coalescent} for $(X^n, \tilde X^n)$:
\[
	X^n_s = \tilde X^n_s
\Quad{implies}
	X^n_t = \tilde X^n_t
\Qforall
	t \ge s.
\]
So, we need to show, for all $\delta>0$ sufficiently small,~that
\[
	\limsup_{n\to\infty}
	\mathbf P_{x^n_0, \tilde x^n_0}\rbr{ \tau_\coal > 4 \delta^{1/2} n }
\le
	2 \delta^{1/4} , 
\] where the initial conditions satisfy \eqref{eq:CloseToBoundary}.
To do so, consider the set
\[
	\Omega'_n
\cq
	\mbz \cap \sbb{ \ka_n^2 n - 2\delta^{-1/4} \sqrt n \ka_n, \ka_n^2 n + 2\delta^{-1/4} \sqrt n \ka_n },
\]
and define $(R_t^n)_{t \ge 0}$ and $(\tilde R_t^n)_{t \ge 0}$ as two birth--death chain on $\Omega'_n$ according to the transition rates $q_n$ from \cref{res:bd:rep}:
ie, for all $x \in \Omega'_n$,
\begin{alignat*}{2}
	q_n(x, x+1)
&=
	2 (k_n - x)^2 / n^2
&&\Quad{provided}
	x+1 \le \ka_n^2 n + 2\delta^{-1/4} \sqrt n \ka_n ,
\\
	q_n(x, x-1)
&=
	2 x (n - 2 k_n + x) / n^2
&&\Quad{provided}
	x-1 \ge \ka_n^2 n - 2\delta^{-1/4} \sqrt n \ka_n 
\end{alignat*}
Let $\sigma_\coal \cq \inf\bra{t \ge 0 \mid R^n_t = \tilde R^n_t}$ and, with a slight abuse of notation, denote by $\tilde{\mathbf P}$ the basic coupling between $(R_t^n)_{t \ge 0}$ and $(\tilde R_t^n)_{t \ge 0}$, ie, we use the same Poisson clocks for both processes whenever they are in the same state, and independent clocks otherwise. In particular, the basic coupling is coalescent for $(R^n, \tilde R^n)$.
We will now relate the times $\sigma_\coal$ and $\tau_\coal$.

We define the exit time $\tau_\partial$ of $\Omega'_n$ for the processes $(X_t^n)_{t \ge 0}$ and $(\tilde{X}_t^n)_{t \ge 0}$ as
\[
	\tau_\partial
\cq
	\inf\bra{t \ge 0 \mid X^n_t \notin \Omega'_n  \vee \tilde X^n_t  \notin \Omega'_n}.
\]
\cref{lem:StochasticDom} yields,
for all $x^n_0 ,\tilde x^n_0 \in [\ka_n^2 n - \delta^{-1/4} \ka_n \sqrt n ,  \ka_n^2 n + \delta^{-1/4} \ka_n \sqrt n ]$,
\[
	\mathbf P_{x^n_0, \tilde x^n_0}
	\rbr{ \tau_\partial > n }
\ge
	1 - \delta^{1/4}.
\label{eq:ExitNeeded}
\]
Let $\mathbf P^*$ be a coupling between the four processes
$(X_t^n)_{t \ge 0}$, $(\tilde{X}_t^n)_{t \ge 0}$, $(R_t^n)_{t \ge 0}$ and $(\tilde R_t^n)_{t \ge 0}$,
where we use the same Poisson clocks between edges in $\Omega'_n$ for all four processes, and the same Poisson clocks for all edges with at least one endpoint outside $\Omega'_n$ for the birth--death chains $(X_t^n)_{t \ge 0}$ and $(\tilde{X}_t^n)_{t \ge 0}$, respectively.
In particular, whenever $x^n_0$ and $\tilde x^n_0$ satisfy \eqref{eq:CloseToBoundary},
\[
	\liminf_{n \to \infty}
	\mathbf P^\ast\rbr{ X^n_s=R^n_s,  \tilde X_s^ n = \tilde R_s^n \: \forall \, s\le n \mid X^n_0 = R^n_0 = x^n_0, \tilde X_s^n = \tilde R_s^n = \tilde x^n_0 }
\ge
	1 - \delta^{1/4}.
\]
Hence, for all $t \le n$,
\[
	\mathbf P_{x^n_0, \tilde x^n_0}\rbr{ \tau_\coal > t }
\le
	\tilde{\mathbf P}_{x^n_0, \tilde x^n_0}\rbr{ \sigma_\coal > t } + \delta^{1/4}.
\]
%
%
Thus, it suffices to show that
for all $\delta>0$ sufficiently small
\begin{equation}
\label{eq:CoalR}
	\limsup_{n\to\infty}
	\tilde{\mathbf P}_{x^n_0, \tilde x^n_0} \big( \sigma_\coal \ge 4 \delta^{1/2} n \big)
\le
	\delta^{1/4}.
\end{equation}

In order to estimate the tails of $\sigma_\coal$, we use a similar strategy to  \cite[Lemma~3.4]{LRS:mc-algs}.
First, we control the drift of the chains $(R_t^n)_{t \ge 0}$ and $(\tilde R_t^n)_{t \ge 0}$ on $\Omega'_n$:
the transition rates~satisfy
\[
	q_n( x,x-1) + q_n(y,y+1) \ge q_n( x,x+1)  1_{\{x+1 \in \Omega'_n\}}  + q_n(y,y-1) 1_{\{y-1 \in \Omega'_n\}}
\Qif
	x > y,
\]
with $x,y \in \Omega'_n$,
so the difference $\abs{ R^n_t - \tilde R^n_t }$ has non-positive drift. 
Hence,
\[
	\Delta^n_t
\cq
	\abs{ R^n_t - \tilde R^n_t }
\Qfor
	t \ge 0
\]
defines a supermartingale with respect to the natural filtration $\mcf$---ie,
\[
	\tilde{\mathbf E}\sbr{ \Delta^n_t \mid \mcf_s } \le \Delta^n_s
\Quad{whenever}
	t \ge s \ge 0 , 
\]
where  $\tilde{\mathbf E}$ denotes the expectation with respect to $\tilde{\mathbf P}$.
We claim that, for $n$ sufficiently large,
\[
	\Phi^n_t
\cq
	(\Delta^n_t)^2 - (4 \delta^{-1/4} \sqrt n \ka_n) \Delta^n_t - \ka_n^2 t
\Qfor
	t \ge 0
\]
defines a submartingle with respect to $\mcf$.
Indeed, the minimal transition rate in $\Omega'_n$ satisfies
\[
	\min_{x \in \Omega'_n}
	\bra{ q_n(x, x-1), q_n(x, x+1) }
\ge
	\ka_n^2
\]
for $n$ sufficiently large, which means that the time-$t$ quadratic variation grows at least~as~$\ka_n^2 t$. Note that $\tilde{\mathbf E}_{x^n_0, \tilde x^n_0}[\tau_\coal^R]< \infty$. 
We apply the Optional Stopping Theorem to the submartingale $(\Phi^n_t)_{t\ge0}$ with stopping time $\tau_\coal^R$ to get for all $x^n_0, \tilde x^n_0$ from \eqref{eq:CloseToBoundary}, with $n$ large enough,  
\begin{equation}\label{eq:ExpCoalR}
\tilde{\mathbf E}_{x^n_0, \tilde x^n_0}[\tau_\coal^R]
\le
	4 \delta^{3/4} n.
\end{equation}
Indeed, $\tilde{\mathbf E}_{x^n_0, \tilde x^n_0}[\Phi^n_\tau] \ge \Phi^n_0$ and the following inequalities hold:
\begin{alignat*}{2}
	\Phi^n_0 &\ge - (\delta \sqrt n \ka_n) (4 \delta^{-1/4} \sqrt n \ka_n) = 4 \delta^{3/4} n \ka_n^2
&&\Qsince
	0 \le \abs{x_0 - \tilde x_0} \le \delta \sqrt n \ka_n;
\\
	\Phi^n_\tau &\le - \ka_n^2 \tau
&&\Qsince
	0 \le \Delta^n_t \le 4 \delta^{-1/4} \sqrt n \ka_n.
\end{alignat*}
We obtain \eqref{eq:CoalR}, and hence \eqref{eq:GoalCoalFinal}, now from \eqref{eq:ExpCoalR} and Markov's inequality.
%
%
\end{Proof}

In the following, we set
\begin{equation}
\label{def:gThetaSigma}
	g(\theta,\sigma)
\cq
	\tv{  N(e^{-\theta/2},\sigma)  -  N(0,1) } .
\end{equation}
Note that the function $(\theta,\sigma) \mapsto g(\theta,\sigma)$ is continuous in both components. This is immediate from the continuity of the density of a Gaussian random variables with respect to its mean and standard deviation. Moreover, recall the Ornstein--Uhlenbeck process $(\mathcal D_t)_{t \ge 0}$ from \cref{lem:ConvergenceSmaller} with initial condition~$z$, and that for any fixed constant $C>0$, we had set
\begin{equation*}
\label{def:t0}
	T^+_n(C)
=
	\tfrac14 n \log n - \tfrac12 n \log C.
\end{equation*}
%

We have now all tools in order to determine the limit profile
when $k \gg \sqrt n$.

\begin{Proof}[Proof of \cref{res:intro:bl} when $k \gg \sqrt n$]
Recall that previously we defined
\[
	\mcs_n
=
	\brb{ x/\rbb{ \sqrt n \ka_n (1-\ka_n) } - \ka_n^2 n \mid x \in \bb{0,k_n} }
\Qand
	Y^n_t
=
	(X^n_{nt} - \ka_n^2 n) / \rbb{ \sqrt n \ka_n (1-\ka_n) }.
\]
Let $C > 0$ and define $\mathcal I_{\eps,C} \cq \sbr{ C(1-\eps), C(1+\eps) }$.
Then, by \cref{res:bd:conc},
\begin{equation}\label{eq:SpaceReduction}
	\limsup_{C \rightarrow \infty}\limsup_{n \to \infty}
	\pr*{ Y^n_{\frac14 \log n - \frac12 \log C} \notin \mathcal I_{\eps,C} \midb X^n_0 = k_n }
=
	0.
\end{equation}
Hence,
in order to determine the limit profile,
it suffices to study the chains $(Y^n_t)_{t \ge 0}$ started from $\mathcal I_{\eps,C}$, shifting the time by $\tfrac14 \log n - \tfrac12 \log C$, for $C$ sufficiently large.

We start with a lower bound on the limit profile. In view of \eqref{eq:SpaceReduction}, it suffices to show that, for all $\eps>0$ and all $\theta\in \mbr$, there exists some $C_0=C_0(\eps,\theta)$ such that, for all $C \ge C_0$,
\begin{equation}\label{eq:LowerBoundGoal}
	\liminf_{n \to \infty}
	\inf_{y \in \mathcal I_{\eps,C} \cap \mathcal S_n}
	\tv{ \pr[y]{Y^n_{\frac12 \log C +\frac14\theta} \in \cdot } - \pr{ Y^n_\infty \in \cdot } }
\ge
	g(\theta,1) - 3\eps,
\end{equation}
where $X^n_\infty$ which is drawn according to the stationary distribution $\pi_n$ of $(X^n_t)_{t\ge0}$, and
\[
	Y^n_\infty
\cq
	\frac{X^n_\infty- \ka^2_n n}{\sqrt n \ka_n(1-\ka_n)}.
\]
In order to establish this,
we take as our distinguishing statistic the set $[\tfrac12 e^{-\theta/2}, \infty)$,
\begin{equation}
\label{eq:LowerBoundStep1}
	\liminf_{n \to \infty}
	\inf_{y \in \mathcal I_{\eps,C} \cap \mcs_n }
	\pr[y]{ Y^n_{\frac12 \log C+ \frac14 \theta} \ge \tfrac12 e^{-\theta/2}}
-	\pr{ Y^n_\infty \ge \tfrac12 e^{-\theta/2} }
\ge
	g(\theta,1) - 3\eps.
\end{equation}
By \cref{lem:StochasticDom}, the infimum in \eqref{eq:LowerBoundStep1} is attained for the initialisation $y^n_0 \cq \min(\mcs_n \cap \mathcal I_{\eps,C})$.
Let $\mathcal D_\infty \sim N(0,1)$ and $\tilde{\mathcal D} \sim N(e^{-\theta/2}, 1)$.
Then, we claim that, for all $C$ large enough,
\begin{align*}
&	\liminf_{n \to \infty}
	\rbb{ 
		\pr*[y^n_0]{ Y^n_{\frac12 \log C + \frac14 \theta} \ge \tfrac12 e^{-\theta/2} }
	-	\pr*{ Y^n_\infty \ge \tfrac12 e^{-\theta/2} }
	}
\\&\qquad
\ge
	\pr*{ \mathcal D_{\frac12 \log C + \frac14 \theta} \ge \tfrac12 e^{-\theta/2} \mid \mathcal D_0 = C(1-\eps) }
-	\pr*{ \mathcal D_\infty \ge \tfrac12 e^{-\theta/2} }
-	2\eps
\\&\qquad
\ge
	\pr{ \tilde{\mathcal D} \ge \tfrac12 e^{-\theta/2} }
-	\pr{ \mathcal D_\infty  \ge \tfrac12 e^{-\theta/2} }
-	3\eps
=
	g(\theta,1) - 3 \eps.
\end{align*}
Indeed, the first inequality follows from \cref{lem:ConvergenceSmaller,cor:StationaryDistribution},
the second from the explicit law of $(\mathcal D_t)_{t \ge 0}$ given in \eqref{eq:OUExplicit} and a continuity argument.
The final equality follows from the definition of $g$ in \eqref{def:gThetaSigma} and the representation of the \TV distance as
\[
	\tv{ \mu - \nu }
=
	\intt{A}
	\rbb{ \tfrac{d\mu}{dx}(x) - \tfrac{d\nu}{dx}(x) }
	dx
=
	\mu(A) - \nu(A)
\Qwhere
	A
\cq
	\brb{ x \midb \tfrac{d\mu}{dx}(x) > \tfrac{d\nu}{dx}(x) };
\]
indeed, if $m > 0$, then the density of $N(m, 1)$ is larger than that of $N(0, 1)$ in $[\tfrac12 m, \infty)$.
Hence \eqref{eq:LowerBoundStep1} is proved, and so the desired lower bound \eqref{eq:LowerBoundGoal} follows.

We now turn to the upper bound on the limit profile.
Let $\eps>0$, and let $(X_t^n)_{t \ge 0}$ and $(\tilde{X}_t^n)_{t \ge 0}$ be two
copies.
We start by showing that there exists a coupling such that
\begin{equation}
\label{eq:UpperBoundStep1}
	\liminf_{n \to \infty}
	\pr*[k_n, \pi_n]{ \abs{ X^n_{\frac14 \log n + \frac14 \theta} - \tilde X^n_{\frac14 \log n + \frac14 \theta} } \le 3 \eps \sqrt n \ka_n }
\ge
	g(\theta,1) - 2 \eps . 
\end{equation}
To do this, let $(\tilde Y^n_t)_{t \ge 0}$ denote the re-scaled version of $(\tilde{X}^n_t)_{t \ge 0}$.
Then, by \cref{cor:StationaryDistribution}, there exists a coupling of $\tilde Y^n_t$ and $\mathcal D_\infty$ such that
\begin{equation}
\label{eq:UpperBoundStep2}
	\liminf_{n \rightarrow \infty}
	\pr*[\tilde X^n_0 \sim \pi_n]{ \tilde Y^n_t \in [\mathcal D_\infty(1-\eps), \mathcal D_\infty(1+\eps)] }
=
	1.
\end{equation}
Similarly, \cref{lem:ConvergenceSmaller,lem:StochasticDom} ensure that there is a coupling such that
\begin{equation}
\label{eq:UpperBoundStep3}
	\inf_{y \in \mathcal S_n \cap \mathcal I_{\eps,C}}
	\pr*[Y^n_0 = y, \mathcal D_0 = C]{ Y^n_{\frac12 \log C + \frac14 \theta} \in [\mathcal D_{\frac12 \log C + \frac14 \theta}(1-\eps), \mathcal D_{\frac12 \log C + \frac14 \theta}(1+\eps)] }
\ge
	1- \eps,
\end{equation}
provided $C$ and $n$ are sufficiently large.
Moreover, there exists a coupling such that
\begin{equation}
\label{eq:UpperBoundStep4}
	\liminf_{C \to \infty}
	\pr*[C]{ \mathcal D_{\frac12 \log C + \frac14 \theta} = \mathcal D_\infty }
\ge
	g(\theta),
\end{equation}
using the explicit description of $(\mathcal D_t)_{t \ge 0}$ from \eqref{eq:OUExplicit} and the coupling representation of the $\TV$-distance.
Officially, each $\eps$, $\theta$ and $C$ requires its own coupling, but this is not an issue.

Recalling \eqref{eq:SpaceReduction}, we combine (\ref{eq:UpperBoundStep2}--\ref{eq:UpperBoundStep4}) to deduce \eqref{eq:UpperBoundStep1}.
We now use \eqref{eq:UpperBoundStep1} to show the remaining upper bound on the limit profile.
Since $\mathcal D_\infty \sim N(0,1)$, we have
\[
	\limsup_{n \to \infty}
	\pr*{ \mathcal D_\infty \in [-\delta^{-1/4},\delta^{-1/4}] }
\le
	2\delta^{1/4}
\Qforall
	\delta > 0.
\]
By \cref{lem:Coalescence} with $\delta=\eps$ sufficiently small, and \eqref{eq:UpperBoundStep1}, there exists a coupling such that
\[
	\limsup_{n \to \infty}
	\pr*[k_n, \pi_n]{ X^n_{\frac14 \log n + \frac14 \theta + 4 \eps^{1/2}} = \tilde X^n_{\frac14 \log n + \frac14 \theta + 4 \eps^{1/2} } }
\le
	g(\theta,1) + 2 \eps + 2\eps^{1/4}.
\]
Since $\eps>0$ was arbitrary, and $\theta \mapsto g(\theta, 1)$ is continuous, we conclude. 
\end{Proof}

\section{Limit Profile for $k \asymp \sqrt n$}
\label{sec:as}

Throughout this section, we assume that $k_n^2 / n \to \alpha \in (0, \infty)$ as $n \to \infty$.
In this regime, the equilibrium distribution $\pi$ is approximately $\Pois(\alpha)$, so has mean and variance order $1$.
There is thus no need to rescale the process to get an `order-$1$~limit'.

\begin{Proof}[Outline of Proof]
We first show that the birth--death chain converges to an $M/M/\infty$ queue in the critical window after rescaling the rates.
Next, we derive an explicit description of the law of an $M/M/\infty$ queue, which we can compare with its $\Pois(\alpha)$ equilibrium distribution.

We calculate the limit profile by applying the concentration result from \cref{res:bd:conc} on the hitting time of te critical window, and then running the $M/M/\infty$ queue.
\end{Proof}

We start by establishing convergence, in a mild sense, to an $M/M/\infty$ queue.
The result can be strengthened to the sense of weak convergence of processes, but we do not need that.

\begin{lem}
\label{res:eq:q-limit}
	%
Suppose that $k_n^2/n \to \alpha \in (0, \infty)$ as $n \to \infty$.
Define
\[
	q_\star(x,y)
\cq
	2 \alpha \one{y = x + 1}
+	2 x \one{y = x - 1 \ge 0}
\Qfor
	x, y \in \mbn.
\]
Suppose that $(x_n)_{n\in\mbn}$ and $(y_n)_{n\in\mbn}$ satisfy $x_n / k_n \to 0$ and $y_n / k_n \to 0$ as $n \to \infty$.
Then,
\[
	n q_n(x_n, y_n)
/
	q_\star(x_n, y_n)
\to
	1
\Qas
	n \to \infty
\Qforall
	x,y \in \mbn.
\]
\end{lem}

\begin{Proof}
This is a straightforward consequence of the birth--death representation \cref{res:bd:rep}:
\[
	n q_n(x_n, x_n + 1)
&=
	2 k_n^2 / n - 4 x_n k_n / n + 2 x_n^2 / n
=
	2 \alpha_n + \oh1,
\\
	n q_n(x_n, x_n - 1)
&=
	2 x_n - 4 x_n k_n / n + 2 x_n^2 / n
=
	2 x_n + \oh1.
\qedhere
\]
\end{Proof}

Now we see that we are, in essence, working with an $M/M/\infty$ queue once we get near the bulk of $\pi$.
As such, our next step is to study the convergence of an $M/M/\infty$ queue.

\begin{lem}
\label{res:eq:q-dist}
Let $\mcx = (\mcx_s)_{s\ge0}$ denote an $M/M/\infty$ queue with arrival rate $\lambda$ and service rate $\mu$.
Suppose that $\mcx_0 = x$.
Then, $\mcx_s \sim B_s + P_s$, where
\[
	B_s
\sim
	\Bin(x, e^{-\mu s})
\Qand
	P_s
\sim
	\Pois\rbb{ \tfrac\lambda\mu (1 - e^{-\mu s}) }
\quad
	\text{independently}.
\]
In particular, the equilibrium distribution is $\Pois(\lambda/\mu)$.
\end{lem}

\begin{Proof}
Customers in an $M/M/\infty$ queue behave independently:
	they are served simultaneously, and have independent $\Exp(\mu)$ service times.
The number $\mcx_s$ in system at time $s$ is precisely
	the number $B_s$ who were there initially and have not left by time $s$
plus
	the number $P_s$ who arrive during $(0,s]$ and have not left by time $s$.
If a customer is in the system at time $r$, then they have probability $\pr{ \Exp(\mu) > s - r } = e^{-\mu(s - r)}$ of still being there at time~$s$.

By independence of the customers' behaviours, $B_s$ and $P_s$ are independent and, moreover, $B_s \sim \Bin(x, e^{-\mu s})$.
Thus, the point process on $[0,s]$ of customers who arrive during $(0,s]$ \emph{and are still there at time $s$} is a Poisson process on $[0,s]$ with \emph{inhomogeneous} rate $\nu(r) \cq \lambda e^{-\mu(s - r)}$.
The number $P_s$ of arrivals in $[0,s]$ of such a Poisson process is Poisson with mean
\[
	\intt[s]{0}
	\nu(r)
	dr
=
	\intt[s]{0}
	\lambda e^{-\mu(s - r)}
	dr
=
	\tfrac\lambda\mu (1 - e^{-\mu s}).
\qedhere
\]
\end{Proof}

In our framework, $\lambda = 2 \alpha$ and $\mu = 2$, so $\lambda/\mu = \alpha$.
We will end up in a scenario where $x = C$ and $e^{-2 s} = e^{-2\theta} / C$, with $C$ large.
In this regime, the Poisson approximation to the Binomial holds:
	$B_s \approx \Pois(e^{-2\theta})$.
The sum of independent Poisson variables is itself Poisson:
\[
	\mcx_s
\approx
	\Pois(e^{-2\theta}) + \Pois\rbb{ \alpha (1 - e^{-2\theta}/C) }
=
	\Pois\rbr{ \alpha + e^{-2\theta} + \alpha e^{-2\theta}/C }
\approx
	\Pois(\alpha + e^{-2\theta}).
\]
This is what gives the limit profile as the \TV distance between $\Pois(\alpha + e^{-2\theta})$ and $\Pois(\alpha)$.
We now make this strategy precise, allowing us to determine the limit profile when $k_n \asymp \sqrt n$.

\begin{Proof}[Proof of Theorem~\ref{res:intro:bl} when $k \asymp \sqrt n$]
We start with a `burn-in' period to get $X^n_{T^-_n(C)} \approx C$:
\begin{equation*}
	\limsup_{C\to\infty}
	\limsup_{n\to\infty}
	\pr{ \abs{ X^n_{T^-_n(C)} - C } > C \eps }
=
	0
\Qforall
	\eps > 0,
\end{equation*}
by \cref{res:bd:conc}, where we recall that
\[
	T^-_n(C) = \tfrac12 n \log k - \tfrac12 n \log C.
\]
Condition on $X^n_{T^-_n(C)} = C (1 + \xi)$, with $\xi \in (-1,1)$, and initialise $\mcx_0 \cq C (1 + \xi)$.
We take
\[
	t_n
\cq
	\tfrac12 n \log k + \theta n,
\Qand
	s
\cq
	\tfrac12  \log C + \theta.
\]
Then, by \cref{res:eq:q-limit},
for all $\theta \in \mbr$,
\[
	\lim_{n \to \infty}
	\tv*{ 
		\pr*{ X^n_{t_n} \in \cdot \mid X^n_{T^-_n(C)} = C (1 + \xi) }
	-	\pr*{ \mcx_s \in \cdot \mid \mcx_0 \cq C (1 + \xi) ) }
	}
=
	0.
\]
Recall that, by \cref{res:eq:q-dist},
\(
	\mcx_s
\sim
	B_s + P_s
\)
where
\[
	B_s
\sim
	\Bin\rbb{ C(1 + \xi), e^{-2s} }
\Qand
	P_s
\sim
	\Pois\rbb{ \alpha (1 - e^{-2s}) }
\quad
	\text{independently.}
\]
Since $e^{-2s} = e^{-2\theta} / C$, the Poisson approximation to the Binomial gives
\[
	\lim_{C \to \infty}
	\tv*{ \pr{ B_s \in \cdot } - \Pois\rbb{ (1 + \xi) e^{-2\theta} } }
=
	0.
\]
Similarly, since $s \to \infty$ as $C \to \infty$, we see that
\[
	\lim_{C \rightarrow \infty}	\tv{ \P( P_s\in \, \cdot \, ) -
	\Pois\rbr{ \alpha }  } = 0  . 
\]
We take $\xi \to 0$ as $C \to \infty$, sufficiently slowly.
Plugging all these estimates in,
\begin{equation*}
	X^n_{t_n}
\stackrel{(n\to\infty)}{\to}
	\mcx_s
\stackrel{(d)}{=}
	\Bin\rbr{ C, e^{-2\theta}/C }
+	\Pois\rbr{ \alpha }
\stackrel{(C\to\infty)}{\to}
	\Pois(\alpha + e^{-2\theta}),
\end{equation*}
with the convergence in \TV.
Recall that here the equilibrium distribution $\pi_n \to \Pois(\alpha)$ in \TV as $n \to \infty$.
Combining these and using the triangle inequality completes the proof.
	%
\end{Proof}

\section{Limit Profile for $k \ll \sqrt n$}
\label{sec:ll}

Throughout this section, we assume that $k_n^2/n \to 0$ as $n \to \infty$.
With help of \cref{res:bd:conc} and the observation that $\pi_n$ is concentrated on $\{0\}$---see, eg, \cite[Equation~(3.1)]{LL:cutoff-ep-ip}---it suffices to give a bound on the time it takes the birth--death chain with rates $q_n$ to hit $0$.

We define for all $x \in \mbn$ the hitting time
\[
	\tau^n_x
\cq
	\inf\bra{ t \geq 0 \mid X^n_t = x }
\]
and, since $\pi_n$ concentrates on $\bra{0}$, we expect that, for all $t \ge 0$, as $n \to \infty$
\[
		\tv{ \pr[k_n]{ X^n_t \in \cdot } - \pi_n } \approx
		\pr*[k_n]{ \tau^n_0 > t }
\]

\begin{Proof}[Outline of Proof]
\qedtriangle
In view of
this intuition,
we proceed in two steps.
\begin{enumerate}[label = \arabic*.]
	\item 
	We argue that the birth--death chain $(X_t^{n})_{t \ge 0}$ started from $x$ can until time $\tau^n_0$ be well-approximated by a sum of independent Exponential random variables.
	
	\item 
	We show that this sum, suitably renormalised, converges to a Gumbel.
\qedhere
\end{enumerate}
\end{Proof}

The following is our main approximation result of this section.

\begin{prop}
\label{pro:HittingLemmaSmall}
	%
For all $\theta \in \mbr$,
\begin{equation}
\label{eq:CoupleHittingTimes}
	\lim_{x \to \infty}
	\lim_{n \to \infty}
	\absb{ 
		\pr[x]{ \tau^n_0 \ge \tfrac12 n \log x + \tfrac12 \theta n }
	-	\pr{ G > \theta }
	}
=
	0,
\end{equation}
where $G \sim \operatorname{Gum}(0,1)$ is a standard Gumbel random variable:
ie,
\[
	\pr{ G \ge x }
=
	1 - e^{-e^{-x}}
\Qforall
	x \in [0, \infty).
\]
\end{prop}

The rates in this regime satisfy $q_n(x,x-1) = 2x/n + \oh{1/n}$ and $q_n(x,x+1) = \oh{1/n}$, where the $\oh1$ term depends on $x$. So, typically, the chain just jumps down until it hits $0$.
Based on this, let $T_z \sim \Exp(z)$ independently, for $z \in \mbn$, and
\[
	S_{x,y}
\cq
	\sumt{z=y+1}[x]
	T_z
\Qfor
	x,y \in \mbn
\Qwith
	y < x.
\]
Then, $\tfrac12 n T_z \sim \Exp(2z/n)$ represents the time to jump $z \to z-1$, and $\tfrac12 n S_{x,y}$ the time~$x \to y$.

\begin{lem}
\label{lem:CouplingExponential}
	%
For all $x,y \in \mbn$ with $y < x$,
\begin{equation}
\label{eq:CoupleDifferences}
	\lim_{n \to \infty}
	\tv{ \pr[x]{ \tau^n_y \in \cdot } - \pr{ \tfrac12 n S_{x,y} \in \cdot } }
=
	0.
\end{equation}
	%
\end{lem}

\begin{Proof}
The chain steps directly $x \to x-1 \to ... \to y+1 \to y$, and the step $z \to z-1$ can be coupled with $\tfrac12 n T_z$, with probability tending to $1$.
The result follows.
\end{Proof}


It is well-known that $S_{m,0}$ has a simple interpretation as a coupon-collection problem.
\begin{quote}
	Consider a collection of $m$ coupons.
	At rate $1$, independently draw coupons, with replacement.
	Then, $S_{m,0}$ has the law of the first time every coupon has been drawn at least once.
\end{quote}
By the memoryless property, this is equivalent to assigning an independent $\Exp(1)$ random variable to each coupon, and taking $S_{m,0}$ to be the maximum of these $m$ coupons.
This allows us to deduce the following classical asymptotic result; see, for example, \cite{ER:coupon-collector}.

\begin{lem}
\label{lem:GumbelApprox}
Let $G \sim \operatorname{Gum}(0,1)$ be a standard Gumbel random variable.
Then,
\[
	\lim_{m \to \infty}
	\rbb{ S_{m,0} - \log m }
\to
	G
\quad
	\text{in distribution}.
\]
\end{lem}

\begin{Proof}
It is easy to see, by the memoryless property, that
\(
	\max_{\ell=1, ..., m}
	Z_\ell
\sim
	S_{m,0}.
\)
Let
\[
	F(t) \cq \pr{\Exp(1) \le t } = 1 - e^{-t}
\Qfor
	t \in [0,\infty),
\]
the cdf of $\Exp(1)$.
Then, for all $x \in \mbr$, by independence,
\[
	\pr{ S_{m,0} \le \log m + x }
&
=
	\pr{ \cap_{\ell=1}^m Z_\ell \le \log m + x }
\\&
=
	F(\log m + x)^m
=
	(1 - e^{-x}/m)^m
\to
	\exp{-e^{-x}}
\Qas
	m \to \infty,
\]
which is the cdf of the standard Gumbel distribution.
\end{Proof}

\begin{Proof}[Proof of Lemma \ref{pro:HittingLemmaSmall}]
This is immediate by combining Lemma \ref{lem:CouplingExponential} and Lemma \ref{lem:GumbelApprox}.
\end{Proof}

We have now all tools to show the remaining part of Theorem~\ref{res:intro:bl}.

\begin{Proof}[Proof of Theorem~\ref{res:intro:bl} when $k \ll \sqrt n$]
Recall from \cref{res:bd:conc} that
\begin{equation}
\label{eq:LocalConvergenceCoupon}
	\limsup_{C\to\infty}
	\limsup_{n\to\infty}
	\pr{ \abs{ X^n_{T_n^-(C)} - C } > C \eps }
=
	0
\Qforall
	\eps > 0,
\end{equation}
where we recall
\[
	T_n^-(C)
=
	\tfrac12 n \log k_n - \tfrac12 \log C.
\]
Since $\pi_n(\bra{0}) = 1 - \oh1$, the following upper and lower bounds hold:
\begin{empheq}
[ left = { \displaystyle
	\tv*{ \pr[k_n]{ X^n_t \in \cdot } - \pi_n }
\empheqlbrace} ]%
{alignat = 1}
	&\le \pr[k_n]{ \tau^n_0 > t } + \oh1,
\label{eq:ll:upper:tv}
\\
	&\ge \pr[k_n]{ X^n_t \ne 0 } - \oh1.
\label{eq:ll:lower:tv}
\end{empheq}
We assume that $C$ is much larger than $\abs \theta$, and evaluate these at
\[
	t
\cq
	\tfrac12 n \log k_n + \tfrac12 \theta n;
\Quad{then,}
	t - T_n^-(C)
=
	\tfrac12 \log C + \tfrac12 \theta n.
\]
Below, we write $\eta_{C,n} = o_{C,n}(1)$ if $\limsup_{C\to\infty} \limsup_{n\to\infty} \eta_{C,n} = 0$.

We start with the upper bound.
Let $\eps > 0$.
By \eqref{eq:LocalConvergenceCoupon} and monotonicity,
\[
	\pr[k_n]{ \tau^n_0 > t }
&
\le
	\pr*[k_n]{ \tau^n_0 > t, \: X^n_{T_n^-(C)} = \ceil{ C(1+\eps) } }
+	o_{C,n}(1)
\\&
\le
	\pr*[\ceil{ C(1+\eps) }]{ \tau^n_0 > t - T_n^-(C) }
+	o_{C,n}(1).
\]
Assume $C \ge 1/\eps$.
Hence, applying \cref{pro:HittingLemmaSmall} with $x \cq \ceil{C(1+\eps)} \le \floor{ C e^{2\eps} }$,
\[
	\pr[k_n]{\tau^n_0 > t}
&
\le
	\pr[\floor{ Ce^{2\eps}}]{ \tau^n_0 > \tfrac12 \log(C e^{2\eps}) + \tfrac12 (\theta - 2 \eps) n }
\\&
\le
	\pr{G > \theta - 2 \eps}
+	o_{C,n}(1).
\label{eq:ll:upper:hitting}
\nt
\]
Taking $n \to \infty$ and then $C \to \infty$,
plugging \eqref{eq:ll:upper:hitting} into \eqref{eq:ll:upper:tv} gives
\[
	\liminf_{n\to\infty}
	\tv*{ \pr[k_n]{ X^n_t \in \cdot } - \pi_n }
\le
	\pr{ G > \theta - 2 \eps }.
\]

We now turn to the lower bound.
Let $\eps \in (0, \tfrac1{10})$.
By \eqref{eq:LocalConvergenceCoupon} and monotonicity again,
\[
	\pr[k_n]{ X^n_t \ne 0 }
&
\ge
	\pr*[k_n]{ X^n_t \ne 0 \midb X^n_{T_n^-(C)} = \floor{ C(1-\eps) } }
+	o_{C,n}(1)
\\&
=
	\pr[\floor{C(1-\eps)}]{ X^n_{t - T_n^-(C)} \ne 0 }
+	o_{C,n}(1)
\\&
\ge
	\pr[\floor{C(1-\eps)}]{ \tau^n_0 > t - T_n^-(C) }
+	o_{C,n}(1).
\]
Assume $C \ge 2/\eps$.
Hence, applying \cref{pro:HittingLemmaSmall} with $x \cq \ceil{C(1-\eps)} \ge \floor{ Ce^{-2\eps}}$,
\[
	\pr[k_n]{ X^n_t \ne 0 }
&
\ge
	\pr[\floor{Ce^{-2\eps}}]{ \tau^n_0 > \tfrac12 \log(C e^{-2\eps}) + \tfrac12 (\theta + 2 \eps) n }
\\&
\ge
	\pr{ G > \theta + 2 \eps }
+	o_{C,n}(1).
\label{eq:ll:lower:hitting}
\nt
\]
Taking $n \to \infty$ and then $C \to \infty$,
plugging \eqref{eq:ll:lower:hitting} into \eqref{eq:ll:lower:tv} gives
\[
	\limsup_{n\to\infty}
	\tv*{ \pr[k_n]{ X^n_t \in \cdot } - \pi_n }
\le
	\pr{ G > \theta + 2 \eps }.
\]

These upper and lower bounds match up to the $\eps$-term:
\[
	\pr{ G > \theta + 2 \eps }
&
\le
	\eqmakebox[liminfsup]{$\displaystyle\liminf_{n\to\infty}$}
	\tv*{ \pr[k_n]{ X^n_t \in \cdot } - \pi_n }
\\&
\le
	\eqmakebox[liminfsup]{$\displaystyle\limsup_{n\to\infty}$}
	\tv*{ \pr[k_n]{ X^n_t \in \cdot } - \pi_n }
\le
	\pr{ G > \theta - 2 \eps }.
\]
But, the function $\theta \mapsto \pr{ G > \theta }$ is continuous, so taking $\eps \to 0$ completes the proof.
\end{Proof}

\section*{Bibliography}
\addcontentsline{toc}{section}{Bibliography}

\renewcommand{\bibfont}{\sffamily\small}
\printbibliography[heading=none]

\end{document}